\documentclass[11pt]{article}

\bibstyle{plain}

\hsize \textwidth      % Longueur d'une ligne
\advance \hsize by -\marginparwidth
\evensidemargin \oddsidemargin % Idem pour les pages de gauche

\advance\hoffset by 5mm % Pour corriger un decalage residuel sur la gauche

\newcommand{\ind}{\textbf{1}}
\usepackage{amssymb}
\usepackage{amsfonts}
\usepackage{amsmath}
\usepackage{amsxtra}
\usepackage[english,german]{babel}
\usepackage{epsfig}
\usepackage{latexsym}

\numberwithin{equation}{section}

\newtheorem{thm}{Theorem}[section]
\newtheorem{prop}[thm]{Proposition}
\newtheorem{cor}[thm]{Corollary}
\newtheorem{lemma}[thm]{Lemma}

\newtheorem{rem}[thm]{Remark}
\newtheorem{example}[thm]{Example}
\newtheorem{Assumption}[thm]{Assumption}

\newenvironment{Proof}{\textsc{Proof:}}{\mbox{ } \hfill $\Box$ \vspace{2mm}}

\newenvironment{rmenumerate}
    {\begin{enumerate}}
    {\end{enumerate}}

\newcommand{\ma}{\mathbb{A}}
\newcommand{\n}{\mathbb{N}}
\newcommand{\re}{\mathbb{R}}
\newcommand{\E}{\mathbb{E}}
\newcommand{\Prob}{\mathbb{P}}

\newcommand{\mf}{\mathcal{F}}

\newcommand{\half}{\frac{1}{2}}

\newcommand{\fbm}{fractional Brownian motion }

 \normalsize

\begin{document}

\selectlanguage{english}
%%%%%%%%%%%%%%%%%%%%%%%%%%%%%%%%%%%%%%%%%%%%%%%%%%%%%%%%
%%%%%%%% Title %%%%%%%%%%%%%%%%%%%%%%%%%%%%%%%%%%%%%%%%%
%%%%%%%%%%%%%%%%%%%%%%%%%%%%%%%%%%%%%%%%%%%%%%%%%%%%%%%%

\title{Queueing Theoretic Approaches to Financial Price
Fluctuations\thanks{We thank W.~Massey and participants of the
workshops on ``Microscopic Stochastic Dynamics in Economics'' and
``Complexity and Randomness in Economic Dynamical Systems'' at
Bielefeld University for valuable comments and discussions.}}

\author{Erhan Bayraktar\thanks{Department of Mathematics,
University of Michigan, 530 Church Street, Ann Arbor, MI 48109,
{\em erhan@umich.edu}.}\and Ulrich Horst \thanks{Department of
Mathematics,  The University of British Columbia,  1984
Mathematics Road Vancouver, BC, V6T 1Z2 Canada; {\em
horst@math.ubc.ca}} \and Ronnie Sircar
\thanks{Operations Research \& Financial Engineering Department, Princeton University, E-Quad, Princeton NJ 08544;
{\em sircar@princeton.edu}.} }

%\date{\today}

\maketitle

\begin{abstract}
One approach to the analysis of stochastic fluctuations in market
prices is to model characteristics of investor behaviour and the
complex interactions between market participants, with the aim of
extracting consequences in the aggregate. This agent-based
viewpoint in finance goes back at least to the work of Garman
(1976) and shares the philosophy of statistical mechanics in the
physical sciences. We discuss recent developments in market
microstructure models. They are capable, often through numerical
simulations, to explain many stylized facts like the emergence of
herding behavior, volatility clustering and fat tailed returns
distributions. They are typically queueing-type models, that is,
models of order flows, in contrast to classical economic
equilibrium theories of utility-maximizing, rational,
``representative'' investors. Mathematically, they are analyzed
using tools of functional central limit theorems, strong
approximations and weak convergence. Our main examples focus on
investor inertia, a trait that is well-documented, among other
behavioral qualities, and modelled using semi-Markov switching
processes. In particular, we show how inertia may lead to the
phenomenon of long-range dependence in stock prices.
\end{abstract}

%\begin{center}
%{\bf To appear in Handbook of Financial Engineering}
%\end{center}

\tableofcontents

%\newpage

\section{Introduction}
Modeling market microstructure in order to understand the effects
of many individual investors on aggregate demand and price
formation is both a classical area of study in economics, and a
rapidly growing activity among researchers from a variety of
disciplines, partly due to modern-day computational power for
large-scale simulations, and the increased availability of price
and order-book data. Among the benefits of this type of analysis,
whether mathematical or simulation-based, is the design of better
models of macroscopic financial variables such as prices, informed
by microscopic (investor-level) features, that can then be
utilized for improved forecasts, investment and policy decisions.

The approach we discuss here is to identify characteristics common
to large groups of investors, for example prolonged inactivity or
inertia, and study the resulting price dynamics created by order
flows. Typically, we are interested in understanding the
microstructure effects on the aggregate quantity through
approximations from stochastic process limit theorems when there
is a large number of investors.

In this point of view, we model right away the behavior of
individual traders rather than characterizing agents' investment
decisions as solutions to individual utility maximization problems.
Such an approach has also been taken in \cite{Garman},
\cite{Foellmer-Schweizer}, \cite{Lux} and \cite{horstkirman}, for
example. As pointed out by O'Hara in her influential book Market
Microstructure Theory \cite{ohara}, it was Garman's 1976 paper
\cite{Garman} that ``inaugurated the explicit study of market
microstructure". There, he explains the philosophy of this approach
as follows: ``The usual development here would be to start with a
theory of individual choice. Such a theory would probably include
the assumption of a stochastic income stream
[and] probabilistic budget constraints $\cdots$. %axioms of rational choice
%amongst stochastic consumption streams, and so forth.
But here we are concerned rather with aggregate market behavior
and shall adopt the attitude of the physicist who cares not
whether his individual particles possess rationality, free will,
blind ignorance or whatever, as long as his statistical mechanics
will accurately describe the behavior of large ensembles of those
particles''. This approach is also common in much of the
econophysics literature (see the discussion in \cite{Farmer-zero},
for example), and is of course prevalent in queueing models of
telephone calls or internet traffic \cite{yao}, where interest is
not so much on causes of phone calls or bandwidth demand, but on
phenomenological models and their overall implications. As one
econophysicist explained it in reaction to the usual battle-cry of
the classical economist about rational behaviour, when AT\&T uses
queuing models, it doesn't ask {\em why} you call your
grandmother.

In this article, we provide an outline to recent surveys on
agent-based computational models and analytical models based on
dynamical systems, while our focus is on developing limit theorems
for queueing models of investor behaviour, which apply modern
methods from stochastic analysis to models based on economic
intuition and empirical evidence. The goal is in obtaining insights
into market dynamics by understanding price formation from typical
behavioral qualities of individual investors.

The remainder of this paper is summarized as follows: in
Section~\ref{sec:micro-structure-model}, we briefly survey some
recent research on agent based models. These models relate the
behavioral qualities of investors and quantitative features of the
stock price process. We give a relevant literature review  of
Queuing Theory approaches to the modeling of stock price dynamics in
Section~\ref{sec:order-book}. In
Section~\ref{sec:inertia-in-fin-markets}, we discuss evidence of
investor inertia in financial markets, and we study its effect on
stock price dynamics  in Section~\ref{Sec-Model}. Key tools are a
functional central limit theorem for semi-Markov processes and
approximation results for integrals with respect to fractional
Brownian motion, that establish a link between investor inertia and
long range dependence in stock price returns. These are extended in
Section~\ref{Sec-Feedback} to allow for the {\em feedback} of price
of the stock into agents' investment decisions, using methods and
techniques from \emph{state dependent queuing networks}. We
establish approximation results for the stock price in a
non-Walrasian framework in which the order rates of the agents
depend on the stock price and exogenously specified investor
sentiment. Section \ref{conc} concludes and discusses future
directions.

%%%%%%%%%%%%%%%%%%%%%%%%%%%%%%%%%%%%%%%%%%%%%%%%%%%%%%%%%%%%%%%%%%%%%%%%

\section{Agent-Based Models of Financial
Markets}\label{sec:micro-structure-model}

In mathematical finance, the dynamics of asset prices are usually
modelled by trajectories of some exogenously specified stochastic
process defined on an underlying probability space $(\Omega, \mf,
\Prob)$. Geometric Brownian motion has long become the canonical
reference model of financial price dynamics. Since prices are
generated by the demand of market participants, it is of interest to
support such an approach by a microeconomic model of interacting
agents. In recent years there has been increasing interest in
agent-based models of financial markets where the demand for a risky
asset comes from many agents with interacting preferences and
expectations. These models are capable of reproducing, often through
simulations, many ``stylized facts'' like the emergence of herding
behavior \cite{horstkirman,Lux}; volatility clustering
\cite{cont-volclustering,LuxMarchesi}, or fat-tailed distributions
of stock returns \cite{CB}, that are observed in financial data.

In contrast to the traditional framework of an economy with a
utility-maximizing representative agent, agent-based models
typically comprise many \textsl{heterogeneous} traders who are
so-called \textsl{boundedly rational}. Behavioral finance models
assume that market participants do not necessarily share identical
expectations about the future evolution of asset prices or
assessments about a stock's fundamental value. Instead, agents are
allowed to use rule of thumb strategies when making their investment
decisions and to switch randomly between them as time passes.
Following the seminal article of Frankel and Froot
\cite{FrankelFroot}, one typically distinguishes
\textsl{fundamentalists}, \textsl{noise traders} and
\textsl{chartists}\footnote{Survey data showing the importance of
chartist trading rules among financial practitioners can be found
in, e.g., \cite{TaylorAllen} and \cite{Frankel-Froot-Data}}. A
fundamentalist bases his forecasts of future asset prices upon
market fundamentals and economic factors such as dividends,
quarterly earnings or GDP growth rates. He invests in assets he
considers undervalued, that is, he invests in assets whose price is
beneath his subjective assessment of the fundamental value.
Chartists, on the other hand, base their trading strategy upon
observed historical price patterns such as trends. Technical traders
try to extrapolate future asset price movements from past
observations. Fundamentalists and chartists typically coexist with
fractions varying over time as agents are allowed to change their
strategies in reaction to either the strategies' performances or the
choices of other market participants. Some of these changes can be
self reinforcing when agents tend to follow the more successful
strategies or agents. This may lead to temporary deviations of
prices from the benchmark fundamental or rational expectations
prices generating bubbles or crashes in periods when technical
trading predominates. Fundamentalists typically have a stabilizing
impact on stock prices.

In this section, we review some agent-based models of financial
markets. Our focus will be on a class probabilistic models in
which asset price dynamics are modelled as stochastic processes in
a random environment of investor sentiment. These models are
perhaps most amenable to rigorous mathematical results. Behavioral
finance models based on deterministic dynamical systems are
covered only briefly as they are discussed extensively in a recent
survey by Hommes \cite{hommes-handbook}. For results on
evolutionary dynamics of financial markets we refer to
\cite{Hens-Schenk}, \cite{EHS}, or \cite{Sandroni} and references
therein.

%%%%%%%%%%%%%%%%%%%%%%%%%%%%%%%%%%%%%%%%%%%%%%%%%%%%%%%%%%%%%%%%%%

\subsection{Stock Prices as Temporary Equilibria in Random
Media}

F\"ollmer and Schweizer \cite{Foellmer-Schweizer} argue that asset
prices should be viewed as a sequence of temporary equilibrium
prices in a random environment of investor sentiment; see also
\cite{Foellmer94}. In reaction to a proposed price $p$ in period
$t$, agent $a \in \ma$ forms a random excess demand
$e^a_t(p,\omega)$, and the actual asset price $P_t(\omega)$ is
determined by the market clearing condition of zero total excess
demand. In \cite{Foellmer-Schweizer}, individual excess demand
involves some \textsl{exogenous} liquidity demand and an
\textsl{endogenous} amount obtained by comparing the proposed
price $p$ with some reference level $\hat{S}^a_t$. This dependence
is linear on a logarithmic scale and individual excess demand
takes the form
\begin{equation}
\label{excess-demand}
    e^a_{t}(p,\omega) := c^a_t(\omega) \left( \log \hat{S}^a_{t}(\omega) - \log p
    \right) + \eta^a_{t}(\omega)
\end{equation}
with nonnegative random coefficients $c^a_t(\omega)$. Here
$\eta^a_t(w)$ is the individual's liquidity demand. The logarithmic
equilibrium price $S_{t}(\omega) := \log P_{t}(\omega)$ is then
determined via the market clearing condition $\sum_{a \in \ma}
e^a(P_t(\omega),\omega) = 0.$ It is thus formed from an aggregate
%\begin{equation}
%\label{average}
%    S_{t}(\omega) = \frac{1}{c_t(\omega)} \sum_{a \in \ma} c^a_t(\omega)
%    \hat{S}^a_{t}(\omega) + \eta_{t}(\omega),
%\end{equation}
of individual price assessments and liquidity demands. If the agents
have no sense of the direction of the market and simply take the
last logarithmic price $S_{t-1}$ as their reference level, i.e., if
$\log \hat{S}^a_t = S_{t-1}$, then the log-price dynamics reduces to
an equation of the form
\[
    S_{t} = S_{t-1} + \eta_{t}
\]
were $\eta_t$ denotes the aggregate liquidity demand. In this case
the dynamics of logarithmic prices reduces to a simple random walk
model if the aggregate liquidity demand is independent and
identically distributed over time. This is just the discretized
version of the Black-Scholes-Samuelson geometric Brownian motion
model.

A \textsl{fundamentalists} bases his investment decision on the
idea that asset prices will move closer to his subjective
benchmark fundamental value $F^a$. In a simple log-linear case
\begin{equation} \label{fundamentalist}
    \log \hat{S}^{a}_{t} := S_{t-1} + \alpha^a_t (F^{a} - S_{t-1})
\end{equation}
for some random coefficient $0 < \alpha^a_t < 1$. If only such
\textsl{information traders} are active on the market, the
resulting logarithmic stock price process takes the form of a
mean-reverting random walk in the random environment
$\{\alpha_t\}_{t \in \n}$ $(\alpha_t = \{\alpha^a_t\}_{a \in
\ma})$.
%\[
%    S_{t} = S_{t-1} + \bar{\alpha}_t (S_{t-1} - \bar{F}_t).
%\]
A combination of information trading and a simple form of
\textsl{noise trading} where some agents take the proposed price
seriously as a signal about the underlying fundamental value
replacing $F^a$ in (\ref{fundamentalist}) by $p$ leads to a class
of discrete time Ornstein-Uhlenbeck processes. Assuming for
simplicity that subjective fundamentals equal zero the logarithmic
price process takes the form
\begin{equation} \label{discrete-time-prices}
    S_{t} - S_{t-1} = \tilde{\gamma}_t S_{t-1} +
    \gamma_t
\end{equation}
with random coefficients $\tilde{\gamma}_t$ and $\gamma_t$. These
coefficients describe the fluctuations in the proportion between
fundamentalist and noise traders. When noise trading predominates,
$\tilde{\gamma}_t$ becomes negative and the price process transient.
Asset prices behave in a stable manner when the majority of the
agents adopts a fundamentalist benchmark.

\subsubsection{Random environment driven by interactive Markov
processes}

Let us now discuss a possible source of randomness driving the
environment for the evolution of stock prices. Extending an
earlier approach in \cite{Foellmer94}, Horst \cite{Horst02}
analyzes a situation with countably many agents located on some
integer lattice $\ma$ where the environment is generated by an
underlying Markov chain with an interactive dynamics. There is a
set $C$ of possible characteristics or trading strategies. An
agent's state $x^a_t \in C$ specifies her reference level for the
following period. The environment is then driven by a Markov chain
\begin{equation} \label{product-kernel}
    \Pi(x_t;\cdot) = \prod_{a \in \ma} \pi_a(x_t;\cdot)
\end{equation}
where $x_t = (x_t^a)_{a \in \ma}$ denotes the current
configuration of reference levels. The distribution
$\pi_a(x_t;\cdot)$ of an agent's state in the following period may
depend both on the current states of some ``neighbors'' and
signals about the aggregate behavior. Information about aggregate
behavior is carried in the empirical distribution $\varrho(x_t)$
or, more generally, the empirical field $R(x_t)$ associated to the
configuration $x_t$. The empirical field is defined as the weak
limit
\[
    R(x_t) := \lim_{n \rightarrow \infty} \frac{1}{|\ma_n|} \sum_{a \in \ma_n} \delta_{\theta^a
    x_t}(\cdot)
\]
along an increasing sequence of finite sub-populations $\ma_n
\uparrow \ma$ and $(\theta^a)_{a \in \ma}$ denotes the canonical
shift group on the space of all configurations. Due to the
dependence of the transition probabilities $\pi_a(x;\cdot)$ on
aggregate behavior, the kernel $\Pi$ does not have the Feller
property, and so standard convergence results for Feller processes
on compact state spaces do not apply. As shown by F\"{o}llmer and
Horst \cite{Foellmer-Horst} and Horst \cite{Horst-AdvApp} the
evolution of aggregate behavior on the level of empirical fields
can be described by a Markov chain. In \cite{Horst02}, it is the
$\{R(x_t)\}_{t \in \n}$ process that generates the environment:
\[
    (\tilde{\gamma}_t,\gamma_t) \sim Z(R(x_t);\cdot) \quad \mbox{
    for a suitable stochastic kernel $Z$.}
\]
Under a weak interaction condition the process $\{R(x_t)\}_{t \in
\n}$ settles down to a unique limiting distribution. Hence asset
prices asymptotically evolve in a stationary and ergodic random
environment. This allows us to approximate the discrete time
process $\{S_t\}_{t \in \n}$ by the unique strong solution to the
stochastic differential equation
\[
    dZ_t = Z_t dX_t + d\tilde{X}_t
\]
where $X$ and $\tilde{X}$ are Brownian motions with drift and
volatility; see \cite{Foellmer-Schweizer} or \cite{Horst02} for
details.

\subsubsection{Feedback effects}

The random environment in \cite{Horst02} is generated by a Markov
process describing the evolution of individual behavior. While
this approach is capable of capturing some interaction and
imitation effects such as word-of-mouth advertising unrelated to
market events, the dynamics of $\{x_t\}_{t \in \n}$ lacks a
dependence on asset price dynamics. The model by F\"{o}llmer,
Horst, and Kirman \cite{horstkirman} captures feedback effects
from stock prices into the environment. At the same time it allows
for trend chasing. A trend chaser or chartist bases his
expectation of future asset prices and hence his trading strategy
upon observed historical price patterns such as trends. In
\cite{horstkirman}, for instance, the chartist's benchmark level
takes the form
\begin{equation} \label{chartist}
    \log \hat{S}^{a}_{t} := S_{t-1} + \beta^a_t (S_{t-1} -
    S_{t-2}).
\end{equation}
A combination of the trading strategies (\ref{fundamentalist}) and
(\ref{chartist}) yields a class of asset price processes that can
be described by a higher order stochastic difference equation. In
\cite{horstkirman}, the agents use one of a number of predictors
which they obtain from financial ``gurus'' to forecast future
price movements. The agents evaluate the gurus' performance over
time. Performances are measured by weighted sums of past profits
the strategies generate. The probability of choosing a given guru
is related to the guru's success. As a result, the configuration
$x_t$ of individual choices at time $t$ is a random variable whose
distribution depends on the current vector of performance levels
$U_{t-1}$. This dependence of the agents' choices on performances
introduces a feedback from past prices into the random
environment. Loosely speaking one obtains a difference equation of
the form (\ref{discrete-time-prices}) where
%It is generated by an
%underlying \textsl{endogenous} process $\{x_t\}_{t \in \n}$
%describing the evolution agents' individual choices $x^a_t$ of
%gurus. %This leads to a dynamics of the form
\[
    (\tilde{\gamma_t},\gamma_t) \sim Z(U_t;\cdot) \quad \mbox{
    for a suitable stochastic kernel $Z$.}
\]
%for the environment.
While prices can temporarily deviate from
fundamental values, the main result in \cite{horstkirman} shows
that the price process has a unique stationary distribution, and
time averages converge to their expected value under the
stationary measure if the impact of trend chasing is weak enough.

\subsubsection{Multiplicity of equilibria}

As argued by Kirman \cite{Kirman92}, in a random economy with many
heterogenous agents, a natural idea of an equilibrium is not a
particular state, but rather a distribution of states reflecting
the proportion of time the economy spends in each of the states.
In the context of microstructure models where liquidity trading or
interaction effects prevent asset prices from converging pathwise
to some steady state, stationary distributions for asset prices
are thus a natural notion of equilibrium. In this sense, the main
result in \cite{horstkirman} may be viewed as an existence and
uniqueness result for equilibria in financial markets with
heterogenous agents. Horst and Wenzelburger
\cite{horst-wenzelburger} study a related model with many small
investors where performances are evaluated according historic
returns or Sharpe ratios. In the limit of an infinite set of
agents the dynamics of asset prices can be described by a path
dependent linear stochastic difference equation of the form
\[
    Y_{t} = A(\varrho_{t-1}) Y_{t-1} + B(\varrho_{t-1},\epsilon_t).
\]
Here $\{\epsilon_t\}_{t \in \n}$ is an exogenous i.i.d sequence of
noise trader demand and $\varrho_{t-1}$ denotes the empirical
distribution of the random vector $Y_0, Y_1, \ldots ,Y_{t-1}$.
While the models shares many of the qualitative features of
\cite{Horst02} and \cite{horstkirman}, it allows for multiple
limiting distributions of asset prices. If the interaction between
different agents is strong enough, asset prices converge in
distribution to a \textsl{random} limiting measure. Randomness in
the limiting distribution may be viewed as a form of market
incompleteness generated by contagious interaction effects.

%%%%%%%%%%%%%%%%%%%%%%%%%%%%%%%%%%%%%%%%%%%%%%%%%%%%%%%%%%%%%%%%%%%%%%

\subsubsection{Interacting agent models in an overlapping generations
framework}

The work in \cite{horst-wenzelburger} is based on earlier work by
B\"{o}hm {\em et al.} \cite{BDW}, B\"{o}hm and Wenzelburger
\cite{BoehmWenzelJEDC}, and Wenzelburger \cite{WenzelJEDC}. These
authors developed a dynamic analysis of endogenous asset price
formations in the context of overlapping generations economies
where agents live for two periods and the demand for the risky
asset comes from young households. They investigate the impact of
different forecasting rules on both asset price and wealth
dynamics under the assumption that agents are myopic and therefore
boundedly rational, mean-variance maximizers. B\"{o}hm {\em et
al.} \cite{BDW} study asset prices and equity premia for a
parameterized class of examples and investigate the role of risk
aversion and of subjective as well as rational beliefs. It is
argued that realistic parameter values explain Mehra and
Prescott's equity premium puzzle (\cite{Mehra-Prescott}). The
model is generalized in \cite{WenzelJEDC} to a model with an
arbitrary number of risky assets and heterogeneous beliefs, thus
generalizing the classical CAPM. A major result is conditions
under which a learning scheme converges to rational expectations
for one investor while other investors have non-rational beliefs.
A second major result is the notion of a \emph{modified market
portfolio} along with a generalization of the security market line
result stating that in a world of heterogeneous myopic investors,
modified market portfolios are \emph{mean-variance efficient in
the classical sense of CAPM}, regardless of the diversity  of
beliefs of other agents. See \cite{BoChi} for a related approach.

%%%%%%%%%%%%%%%%%%%%%%%%%%%%%%%%%%%%%%%%%%%%%%%%%%%%%%%%%%%%%%%%%%%%%

\subsubsection{Feedback Effects from Program Trading, Large Agents and Illiquidity}
A different type of feedback effect, from the actions of a large
group of {\em program traders} or large influential agents has
been modelled in the financial mathematics literature. In the
1990s, following the Brady report that attributed part of the
cause of the 1987 stock market crash to program trading by
institutions following portfolio insurance strategies, researchers
analyzed the feedback effect from option Delta-hedging by a
significant fraction of market participants on the price dynamics
of the underlying security. See, for example, Frey and Stremme
\cite{frey}, Sircar and Papanicolaou \cite{feedback},
Sch\"{o}nbucher and Wilmott \cite{schonbucher} and Platen and
Schweizer \cite{schweizer}.

Related analyses can be found in models where there is a large
investor whose actions move the price, for example Jonsson and
Keppo \cite{jonsson-keppo}, and where there is a market depth
function describing the impact of order size on price, for example
Cetin {\em et al.} \cite{cetin}. A cautionary note on all such
models is that, under sensible conditions, they do not explain the
implied volatility smile/skew that is observed in modern options
markets (in fact they predict a reverse smile). This would suggest
that program trading, large agent or illiquidity effects are
second order phenomena as far as derivatives markets are
concerned, compared with the impacts of jumps or stochastic
volatility.

There has also been some recent empirical work on estimating the
market depth function, in particular the tail of the distribution
governing how order size impacts trading price: see Farmer and
Lillo \cite{farmer-QF} and Gabaix {\em et al.} \cite{gabaix}.

%%%%%%%%%%%%%%%%%%%%%%%%%%%%%%%%%%%%%%%%%%%%%%%%%%%%%%%%%%%%%%%%%%%%%

\subsection{Stock Prices and Random Dynamical Systems}

An important branch of the literature on agent-based financial
market models analyzes financial markets in which the dynamics of
asset prices can be described by a deterministic dynamical system.
The idea is to view agent-based models as highly nonlinear
deterministic dynamical systems and markets as \textsl{complex
adaptive systems}, with the evolution of expectations and trading
strategies coupled to market dynamics. Many such models, when
simulated, generate time paths of prices which switch from one
expectations regime to another generating \textsl{rational routes to
randomness}, i.e., chaotic price fluctuations. As these models are
considerably more complex than the ones reviewed in the previous
section, analytical characterizations of asset price processes are
typically not available. However, when simulated, these model
generate much more realistic time paths of prices explaining many of
the stylized facts observed in real financial markets.

Particularly relevant contributions include the early work of Day
and Huang \cite{DayHuang}, Frankel and Froot \cite{FrankelFroot} and
the work of Brock and Hommes \cite{BrockHommes}. The latter studies
a model in which boundedly rational agents can use one of two
forecasting rules or investment strategies. One of them is costly
but when all agents use it, the emerging price process is stable.
The other is cheaper but when used by many individuals induces
unstable behavior of the price process. Their model has periods of
stability interspersed with bubble like behavior. In
\cite{BrockHommes98} the same authors introduced the notion of
\textsl{Adaptive Belief Systems} (ABS), a ``financial market
application of the evolutionary selection of expectation rules''
analyzed in \cite{BrockHommes}. An ABS may be viewed as asset
pricing models derived form mean-variance optimization with
heterogenous beliefs. As pointed out in \cite{hommes-handbook}, ``a
convenient feature of an ABS is that it can be formulated in terms
of (price) deviations from a benchmark fundamental and (...) can
therefore be used in experimental and empirical testing of
deviations from the (rational expectations) benchmark.'' Recently,
several modifications of ABSs have been studied. While in
\cite{BrockHommes98} the demand for a risky asset comes from agents
with constant absolute risk aversion utility functions and the
number of trader types is small, Chiarella and He \cite{CH01} and
Brock, Hommes, and Wagener \cite{BHW2005} developed models of
interaction of portfolio decisions and wealth dynamics with
heterogeneous agents whose preferences are described by logarithmic
CRRA utility functions and many types of traders, respectively.
Gaunersdorfer \cite{Gaunersdorfer} extends the work in
\cite{BrockHommes} to the case of time-varying expectations about
variances of conditional stock returns\footnote{There are many other
papers utilizing dynamical system theory to analyze asset price
dynamics in behavioral finance models. For a detailed survey, we
refer the interested reader to \cite{hommes-handbook}.}.

%%%%%%%%%%%%%%%%%%%%%%%%%%%%%%%%%%%%%%%%%%%%%%%%%%%%%%%%%%%%%%%%%%%%%%

\subsection{Queuing Models and Order Book
Dynamics}\label{sec:order-book}

The aforementioned models differ considerably in their degree of
complexity and analytical tractability, but they are all based on
the idea that asset price fluctuations can be described by a
sequence of temporary price equilibria. All agents submit their
demand schedule to a market maker who matches individual demands in
such a way that markets clear. While such an approach is consistent
with dynamic microeconomic theory, it should only be viewed as a
first steps towards a more realistic modelling of asset price
formation in large financial markets. In real markets, buying and
selling orders arrive at different points in time, and so the
economic paradigm that a Walrasian auctioneer can set prices such
that the markets clear at the end of each trading period typically
does not apply. In fact, almost all automated financial trading
systems function as continuous double auctions. They are based on
electronic \textsl{order books} in which all unexecuted limit orders
are stored and displayed while awaiting execution. While
analytically tractable models of order book dynamics would be of
considerable value, their development has been hindered by the
inherent complexity of limit order markets. So far, rigorous
mathematical results have only been established under rather
restrictive assumptions on aggregate order flows by, e.g., Mendelson
\cite{Mendelson}, Luckock \cite{Luckock} and Kruk \cite{Kruk}.
Statistical properties of continuous double auctions are often
analyzed in the econophysics literature e.g., Smith {\em et al.}
\cite{Farmer} and references therein.

Microstructure models with asynchronous order arrivals where
orders are executed immediately rather than awaiting the arrival
of a matching order and where asset prices move into the order to
market imbalance are studied by, e.g. Garman \cite{Garman}; Lux
\cite{LuxEJ,Lux,LuxJEDC} or Bayraktar {\em et al.}
(\cite{kn:bhs}). These models may be viewed as an intermediate
step towards a more realistic modeling of electronic trading
systems.

A convenient mathematical framework for such models, which we will
develop in detail in Section \ref{Sec-Feedback}, is based on the
theory of state-dependent queuing networks (see \cite{MMR} or
\cite{Mandelbaum-Pats} for detailed discussions of Markovian
queuing networks). Underlying this approach is the idea that the
dynamics of order arrivals follows a Poisson-type process with
price dependent rates and that a buying (selling) order increases
(decreases) the stock price by a fixed amount (on a possibly
logarithmic scale to avoid negative prices).

More precisely, the arrival times of aggregate buying and selling
orders are specified by independent Poisson processes $\Pi_+$ and
$\Pi_-$ with price and time dependent rates $\lambda_+$ and
$\lambda_-$, respectively, that may also depend on investor
characteristics or random economic fundamentals. In the simplest
case the logarithmic price process $\{S_t\}_{t \geq 0}$ takes the
form
\[
    S_t = S_0 + \Pi_+ \left( \int_0^t \lambda_+ (S_u,u) du \right)
    - \Pi_- \left( \int_0^t \lambda_- (S_u,u) du \right).
\]
The excess order rate $\lambda_+ (S_u,u) - \lambda_- (S_u,u)$ may be
viewed as a measure of aggregate excess demand while $\Pi_+ \left(
\int_0^t \lambda_+ (S_u,u) du \right) - \Pi_- \left( \int_0^t
\lambda_- (S_u,u) du \right)$ denotes the accumulated net order flow
up to time $t$. In a model with many agents and after suitable
rescaling the asset price process may be approximated by a
deterministic process while the fluctuations around this first order
approximation can typically be described by an Ornstein-Uhlenbeck
diffusion.

Recently, such queuing models have also been applied to modeling
the credit risk of large portfolios by Davis and Esparragoza
\cite{Davis}. They approximate evolution of the loss distribution
of a large portfolio of credit instruments over time. We further
elaborate on queuing theoretic approaches to stock price dynamics
in Section \ref{Sec-Model}. Before that, we introduce a common
investor trait, investor inertia, and show the effects of this
common trait on stock prices.

%%%%%%%%%%%%%%%%%%%%%%%%%%%%%%%%%%%%%%%%%%%%%%%%%%%%%%%%%%%%%%%%%%%

\subsection{Inertia in Financial
Markets}\label{sec:inertia-in-fin-markets}

The models mentioned previously assume that agents trade the asset
in each period. At the end of each trading interval, agents update
their expectations for the future evolution of the stock price and
formulate their excess demand for the following period. However,
small investors are not so efficient in their investment
decisions: they are typically inactive and actually trade only
occasionally. This may be because they are waiting to accumulate
sufficient capital to make further stock purchases; or they tend
to monitor their portfolios infrequently; or they are simply
scared of choosing the wrong investments; or they feel that as
long-term investors, they can defer action; or they put off the
time-consuming research necessary to make informed portfolio
choices. Long uninterrupted periods of inactivity may be viewed as
a form of investor inertia.

\subsubsection{Evidence of inertia}

Investor inertia is a common experience and is well documented.
The New York Stock Exchange (NYSE)'s survey of individual
shareownership in the United States, ``Shareownership2000''
\cite{NYSE}, demonstrates that many investors have very low levels
of trading activity. For example they find that ``23 percent of
stockholders with brokerage accounts report no trading at all,
while 35 percent report trading only once or twice in the last
year''. The NYSE survey also reports (Table 28) that the average
holding period for stocks is long, for example 2.9 years in the
early 90's. Empirical evidence of inertia also appears in the
economic literature. For example, Madrian and Shea \cite{madrian}
looked at the reallocation of assets in employees' individual
401(k) (retirement) plans and found ``a status quo bias resulting
from employee procrastination in making or implementing an optimal
savings decision.'' A related study by Hewitt Associates (a
management consulting firm) found that in 2001, four out of five
plan participants did not do any trading in their 401(k)s. Madrian
and Shea explain that ``if the cost of gathering and evaluating
the information needed to make a 401(k) savings decision exceeds
the short-run benefit from doing so, individuals will
procrastinate.'' The prediction of Prospect Theory (see
\cite{kahneman}) that investors tend to hold onto losing stocks
too long has also been observed in \cite{shefrin}. Another typical
cause is that small investors seem to find it difficult to reverse
investment decisions, as is discussed even in the popular press. A
recent newspaper column (by Russ Wiles in the Arizona Republic,
November 30, 2003) states: ``Perhaps more than anything, investor
inertia is a key force (in financial markets). When the news turns
sour, people tend to hold off on buying rather than bail out. In
2002, the toughest market climate in a generation and a year with
ample Wall Street scandals, equity funds suffered cash outflows of
just one percent.''

\subsubsection{Inertia and long range dependencies in financial time series}
One of the outcomes of a limit analysis of an agent-based model of
investor inertia is a stock price process based on fractional
Brownian motion, which exhibits long-range dependence (that is
correlation or memory in returns). This is discussed in Section
\ref{sec:without-feedback}. In particular, the limit fluctuation
process is a {\em fractional Brownian motion}.

We recall that fractional Brownian motion $B^H$ with
\textsl{Hurst} parameter $H \in (0,1]$ is an almost surely
continuous and centered Gaussian process with auto-correlation
\begin{equation}\label{eq:autocor}
    \mathbb{E}\left\{B^{H}_{t}B^{H}_{s}\right\} =
    \frac{1}{2}\left(|t|^{2H}+|s|^{2H}-|t-s|^{2H}\right).
\end{equation}
\begin{rem}
    Note that the case $H=\frac{1}{2}$ gives standard Brownian
    motion. Also note that the auto-correlation function is
    positive definite if and only if $H \in (0,1]$.
\end{rem}

\begin{figure}[htb]
  \centering
  \includegraphics[totalheight=3in, width=4in]{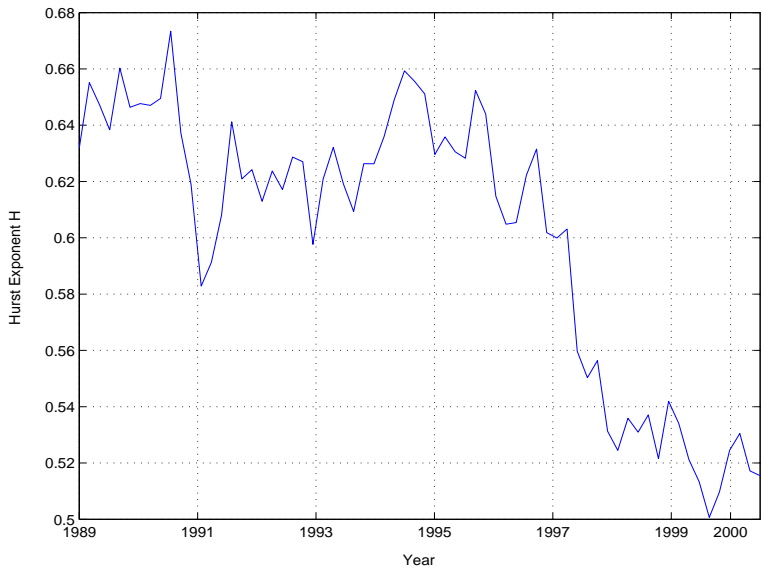}
  \caption{\small{\em Estimates of the Hurst exponent of the S\&P 500 index
  over 1990s, taken from Bayraktar, Poor and Sircar \cite{bps}}.}
  \label{Hests}
\end{figure}

Bayraktar {\em et al.} \cite{bps} studied an asymptotically
efficient wavelet-based estimator for the Hurst parameter, and
analyzed high frequency S\&P 500 index data over the span of 11.5
years (1989-2000). It was observed that, although the Hurst
parameter was significantly higher than the efficient markets
value of $H=\frac{1}{2}$ up through the mid-1990s, it started to
fall to that level over the period 1997-2000 (see Figure
\ref{Hests}).
This might be explained by the increase in Internet trading in
that period, which is documented, for example, in NYSE's
``Shareownership2000'' \cite{NYSE}, Barber and Odean
\cite{odean1}, and Choi {\em et al.} \cite{choi}, in which it is
demonstrated that ``after 18 months of access, the Web effect is
very large: trading frequency doubles.'' Indeed, as reported in
\cite{odean2}, ``after going online, investors trade more
actively, more speculatively and less profitably than before".
Similar empirical findings to that of \cite{bps} were recently
reached, using a completely different statistical technique by
Bianchi \cite{bianchi}.

Thus, the dramatic fall in the estimated Hurst parameter in the
late 1990s can be thought of as {\em a posteriori} validation of
the link the limit theorem in \cite{kn:bhs} provides between
investor inertia and long-range dependence in stock prices. We
review this model in Section \ref{sec:without-feedback}. An
extension based on state dependent queuing networks with
semi-Markov switching is discussed in Section \ref{Sec-Feedback}.

%%%%%%%%%%%%%%%%%%%%%%%%%%%%%%%%%%%%%%%%%%%%%%%%%%%%%%%%%%%%%%%%%%%%%%%%

\section{Microstructure Models with Inert Investors} \label{Sec-Model}
We illustrate the use of microstrucure, or agent-based models,
combined with limit theorems by focusing on investor inertia as a
very common characteristic among small and casual market
participants. In Section \ref{sec:without-feedback} we summarize
earlier work \cite{kn:bhs} that established a mathematical link
between inertia, long-range dependence in stock returns and
potential short-lived arbitrage opportunities for other
`sophisticated' parties. Section \ref{Sec-Feedback} contains an
extension allowing for feedback effects from current prices into
the agents' order rates.

\subsection{A Microstructure Model without Feedback}\label{sec:without-feedback}

We now introduce the basic concepts and notation of the market
microstructure model analyzed in \cite{kn:bhs} that will serve as
basis for the more sophisticated model in
Section~\ref{Sec-Feedback}. We start with a financial market with
a set $\ma:=\{a_1,a_2, \ldots ,a_N\}$ of \textsl{agents} trading a
single risky asset. Each agent $a \in \ma$ is associated with a
continuous-time stochastic process $x^a = \{x^a_t\}_{t \geq 0}$ on
a finite state space $E$ describing his \textsl{trading activity}.

We take a pragmatic approach to specify the demand. Instead of
formulating an individual optimization problem under budget
constraints for the agents, we start right away with the agent's
order rates. The agent $a \in \ma$ accumulates the asset at a rate
$\Psi_t x^a_t$ at time $t \geq 0$. Here $x^a_t$ may be negative
indicating that the agent is selling. The random process $\Psi =
\{\Psi_t\}_{t \geq 0}$ describes the evolution of the size of a
typical trade. It can also be interpreted as a stochastic
elasticity coefficient (the reaction of the price to the market
imbalance). We assume that $\Psi$ is a continuous non-negative
semi-martingale which is independent of the processes $x^a$ and
that $0 \in E$. The agents
do not trade at times when $x^a_t = 0$. % Hence $0$ will be referred
%to as the \textsl{inactive state}.
The holdings of the agent $a \in \ma$ and the ``market imbalance''
at time $t \geq 0$ are thus given by, respectively,
\begin{equation}
\label{process-hat-X}
    \int_0^t \Psi_s x^a_s ds \qquad \mbox{and} \qquad
    \sum_{a \in \ma} \int_0^t \Psi_s x^a_s ds.
\end{equation}
%Hence the process $(I^N_t)_{t \geq 0}$ describes the
%stochastic evolution of the \textsl{market imbalance}.

\begin{rem}
In our continuous time model, buyers and sellers arrive at
different points in time. Hence the economic paradigm that a
Walrasian auctioneer can set prices such that the markets clear at
the end of each trading period does not apply. Rather, temporary
imbalances between demand and supply will occur. Prices will
reflect the extent of market imbalance.
\end{rem}

All the orders are received by a single market maker. The market
maker clears all trades and prices in reaction to the evolution of
market imbalances, the only component driving asset prices.
Reflecting the idea that an individual agent has diminishing
impact on market dynamics if the number of traders is large, we
assume that the impact of an individual order is inversely
proportional to the number of possible traders: a buying (selling)
order increases (decreases) the logarithmic stock price by $1/N$.
The pricing rule for the evolution of the logarithmic stock price
process $S^N=\{S_t^N\}_{t \geq 0}$ is linear and taken to be:
\begin{equation}
\label{process-S}
    dS^N_t = \frac{1}{N} \sum_{a \in \ma} \Psi_t x^a_t dt. %\quad \mbox{and so} \quad
%    S^N_t = S_0 + I_t^N.
\end{equation}

In order to incorporate the idea of market inertia, the agents'
trading activity is modelled by independent and identically
distributed \textsl{semi-Markov} processes $x^a$. Semi-Markov
processes are tailor-made to model individual traders' inertia as
they generalize Markov processes by removing the requirement of
exponentially distributed, and therefore thin-tailed, holding (or
sojourn) times. Since the processes $x^a$ are independent and
identically distributed, it is enough to specify the dynamics of
some ``representative'' process $x=\{x\}_{t \geq 0}$.

%%%%%%%%%%%%%%%%%%%%%%%%%%%%%%%%%%%%%%%%%%%%%%%%%%%%%%%%%%%%%%%%%

\subsubsection{Semi-Markov Processes}

A semi-Markov process $x$ defined on a probability space
$(\Omega,\mf,\Prob)$ is specified in terms of random variables
$\xi_n: \Omega \rightarrow E$ and $T_n: \Omega \rightarrow \re_+$,
satisfying $0 = T_1 \leq T_1 \leq \cdots$ almost surely and
\[
    \Prob \{\xi_{n+1} = j, T_{n+1}-T_{n} \leq t \big|
    \xi_{1},...,\xi_{n}; T_{1},...,T_{n} \} =
    \Prob \{\xi_{n+1} = j, T_{n+1}-T_{n}\leq t \big|\xi_{n}\}
\]
for each $n \in \n$, $j \in E$ and all $t \in \re_+$, through the
relation
\begin{equation}\label{semi-m}
    x_t = \sum_{n \geq 0} \xi_n %\textbf{1}
    \ind_{[T_n,T_{n+1})}(t).
\end{equation}
In economic terms, the representative agent's mood in the random
time interval $[T_{n}, T_{n+1})$ is given by $\xi_n$. The
distribution of the length of the interval $T_{n+1} - T_n$ may
depend on the sequence $\{\xi_n\}_{n \in \n}$ through the states
$\xi_{n}$ and $\xi_{n+1}$. This allows us to assume different
distributions for the lengths of the agents' active and inactive
periods, and in particular to model inertia as a heavy-tailed
sojourn time in the zero state.

\begin{rem}\label{twofour}
In the present analysis of investor inertia, we do not allow for
feedback effects of prices into agents' investment decisions.
While such an assumption might be justified for small,
non-professional investors, it is clearly desirable to allow
active traders' investment decisions to be influenced by asset
prices. We discuss such an extension in the next section.
\end{rem}

We assume that $x$ is temporally homogeneous under the measure
$\Prob$, that is,
\begin{equation}
\label{time-homogeneous}
    Q(i,j,t) \triangleq \Prob \{\xi_{n+1} = j, ~T_{n+1}-T_{n} \leq t \big| \xi_{n} =
    i\}
\end{equation}
is independent of $n \in \n$. By Proposition 1.6 in
\cite{Cinlar75}, this implies that $\{\xi_n\}_{n \in \n}$ is a
homogeneous Markov chain on $E$ whose transition probability
matrix $(p_{ij})$ is given by
\[
    p_{ij} = \lim_{t \rightarrow \infty}Q(i,j,t).
\]
Clearly, $x$ is an ordinary temporally homogeneous Markov process
if $Q$ takes the form
\begin{equation} \label{Markov-kernel}
    Q(i,j,t) = p_{ij} \left(1 - e^{-\lambda_i t} \right).
\end{equation}
We also assume that the \textsl{embedded Markov chain}
$\{\xi_n\}_{n \in \n}$ satisfies $p_{ij} > 0$ so that
$\{\xi_n\}_{n \in \n}$ has a unique stationary distribution. The
conditional distribution function of the length of the $n$-th
sojourn time, $T_{n+1} - T_n$, given $\xi_{n+1}$ and $\xi_n$ is
specified in terms of the \textsl{semi-Markov kernel} $\{
Q(i,j,t); i,j \in E, ~ t \geq 0\}$ and the transition matrix $P$
by
\begin{equation} \label{def-G}
    G(i,j,t) := \frac{Q(i,j,t)}{p_{ij}} = \Prob\{T_{n+1} - T_n \leq t |
    \xi_{n} = i, ~ \xi_{n+1} = j\}.
\end{equation}

The semi-Markov processes are assumed to satisfy the following
conditions.

%The following assumption reflects the idea of market inertia: the
%probability of long uninterrupted trading periods is small compared
%to the probability of an individual agent being inactive for a long
%time.

\begin{Assumption} \label{Assumption-tails}
\begin{rmenumerate}
    \item The average sojourn time at state $i \in E$ is finite:
\begin{equation}
\label{mean-time}
    m_i := \E[T_{n+1} - T_n | \xi_n=i] < \infty.
\end{equation}
    Here $\E$ denotes the expectation operator with respect to
    $\Prob$.
    \item There exists a constant $1 < \alpha < 2$ and a locally bounded function $L: \re_+
    \rightarrow \re_+$ which is slowly varying at infinity (e.g. $\log$), i.e.,
\[
    \lim_{t \rightarrow \infty} \frac{L(xt)}{L(t)} = 1 \qquad \mbox{for
    all} \qquad x>0,
\]
     such that
\begin{equation}
\label{fat-tail}
    \Prob \{ T_{n+1} - T_{n} \geq t \big| \xi_{n} = 0 \}
    \sim t^{-\alpha}L(t).
\end{equation}
 Here we use to notation $f(t) \sim g(t)$ for two functions $f,g:
\re_+ \rightarrow \re_+$ to mean that $\lim_{t\rightarrow
\infty}f(t)/g(t)=1$.

    \item The distributions of the sojourn times at state $i \neq 0$ satisfy
\begin{equation} \label{thin-tail}
    \lim_{t \rightarrow 0} \frac{\Prob\{T_{n+1} - T_{n} \geq t \big| \xi_{n} = i\}}
    {t^{-(\alpha+1)} L(t)} = 0.
\end{equation}
    \item The distribution of the sojourn times in the various
    states have continuous and bounded densities with respect to Lebesgue measure
    on $\re_+$.
\end{rmenumerate}
\end{Assumption}

The key parameter is the tail index $\alpha$ of the sojourn time
distribution of the inactive state zero. Condition
(\ref{fat-tail}) is satisfied if, for instance, the length of the
sojourn time at state $0 \in E$ is distributed according to a
Pareto distribution. The idea of inertia is then reflected by
(\ref{thin-tail}): the probability of long uninterrupted trading
periods is small compared to the probability of an individual
agent being inactive for a long time. In fact, it is natural to
think of the sojourn times in the various active states as being
thin tailed as in the exponential distribution since small
investors typically do not trade persistently.

%%%%%%%%%%%%%%%%%%%%%%%%%%%%%%%%%%%%%%%%%%%%%%%%%%%%%%%%%%%%%%%%%%%%%%%%%%%

\subsubsection{A Limit Theorem for Financial Markets with Inert
Investors}\label{sec:investor-inertia}

We assume that the semi-Markov processes $x^a$ are stationary.
Stationarity can be achieved by a suitable specification of the
common distribution of the initial states and initial sojourn
times. We denote the resulting measure on the canonical path space
by $\Prob^*$. Independence and stationarity of the semi-Markov
processes guarantees that the logarithmic price process can be
approximated pathwise by the process $\{s_t\}_{t \geq 0}$ defined
by
\[
    s_t = \mu \int_0^t \Psi_s ds \quad \mbox{where} \quad \mu :=
    \E^* x^a_0
\]
when the number of agents grows to infinity. Our functional
central limit theorem for stationary semi-Markov processes shows
that after suitable scaling, the fluctuations around $(s_t)_{t
\geq 0}$ can be approximated in law by a process with long range
dependence. The convergence concept we shall use is weak
convergence with respect to the measure $\Prob^*$ of the Skorohod
space $\mathbb{D}$ of all right continuous processes. We write
$\mathcal{L}$-$\lim_{n\rightarrow \infty} Y^n = Y$ if $\{Y^n\}_{n
\in \n}$ is a sequence of $\mathbb{D}$-valued stochastic processes
that converges weakly to the process $Y$.

The convergence result is formulated in terms of a scaling limit
for the processes $\{x^a_{Tt}\}_{t \geq 0}$ $(T \in \n)$. For $T$
large, $x^a_{Tt}$ is a ``speeded-up" semi-Markov process. In other
words, the investors' individual trading dispensations are
evolving on a faster scale than $\Psi$. Observe, however, that we
are not altering the main qualitative feature of the model: agents
still remain in the inactive state for relatively much longer
times than in an active state. In the rescaled model the
logarithmic asset price process $S^{N,T}$ is given by %aggregate order rate
%at time $t$ is given by
\begin{equation} %\label{rescaled-price}
    S^{N,T}_t = \frac{1}{N} \int_0^t \sum_{a \in \ma} \Psi_u x^a_{Tu} du \label{YepsN}.
\end{equation}
The central limit theorem allows us to approximate the
fluctuations around the first order approximation as $N
\rightarrow \infty$. In terms of the Gaussian processes $X^T$ and
$Y^T$ defined by
\begin{equation} \label{def-Gauss-processes}
    X^{T}_{t} \triangleq \mathcal{L}\mbox{-}\lim_{N \rightarrow \infty}
    T^{1-H} \frac{1}{\sqrt{N}}
    \sum_{a=1}^N (x^a_{Tt} - \mu t) \quad \mbox{and} \quad Y^T_t \triangleq
    \int_0^t X^T_s ds,
\end{equation}
with $H=(3-\alpha)/2$, the fluctuations around the first order
approximation can be approximated by an integral of the elasticity
coefficient with respect to $Y^T$:
\[
    \mathcal{L}\mbox{-}\lim_{N \rightarrow \infty}
    \sqrt{N} \left\{ S^{N,T}_t - \mu t \right\}_{0 \leq t \leq 1}
    =  \left\{ \int_0^t \Psi_s d Y^T_s \right\}_{0 \leq t \leq 1}.
\]
In order to see more clearly the effects of investor inertia, we
rescale the price process in space and time and $T$ tend to
infinity. In a benchmark model with many agents where $\Psi \equiv
1$ these, fluctuations when suitably normalized, can be described
by a fractional Brownian motion $B^H$ if $T \rightarrow \infty$.
The Hurst coefficient is related to the degree of investor
inertia.

\begin{thm} \label{thm2} (\cite{kn:bhs})
Let $H = \frac{3-\alpha}{2}$. Assume that $\Psi \equiv 1$, that
Assumption~\ref{Assumption-tails} holds and that $\mu \neq 0$.
Then there exists $\sigma > 0$ such that
\begin{equation} \label{limit1}
    \mathcal{L}\mbox{-}\lim_{T \rightarrow \infty} \mathcal{L}\mbox{-}\lim_{N
    \rightarrow \infty} T^{1-H}\frac{\sqrt{N}}{\sqrt{L(T)}}
    \left\{ S^{N,T}_t - \mu t \right\}_{0 \leq t \leq 1} =  \left\{
    \sigma B^H_t \right\}_{0 \leq t \leq 1}
\end{equation}
\end{thm}

To generalize this result to a market in which the agents' order
rates are coupled by a stochastic elasticity coefficient as in
(\ref{process-S}), we need the following approximation result for
stochastic integrals of continuous semi-martingales with respect
to fractional Brownian motion.

\begin{thm} \label{thm3} (\cite{kn:bhs})
Let $\{\Psi^n\}_{n \in \n}$ be a sequence of good semimartingales
and $\{Z^n\}_{n \in \n}$ be a sequence of $\mathbb{D}$-valued
stochastic processes that satisfy
\begin{itemize}
    \item[(i)] The sample paths of the processes $Z^n$ are almost surely of
    zero quadratic variation on compact sets, and $\Prob \{Z^n_0 = 0\}
    = 1$.
    \item[(ii)] The stochastic integrals $\int \Psi^n d Z^n$ and $\int
    Z^n dZ^n$ exist as limits in probability of Stieltjes-sums, and
    the sample paths $t \mapsto \int_0^t Z^n_s d Z^n_s$ and $t \mapsto
    \int_0^t \Psi^n_s d Z^n_s$ are c\`{a}dl\`{a}g.
\end{itemize}
If $\Psi$ is a continuous semimartingale and if $B^H$ is a \fbm
process with Hurst parameter $H > \frac{1}{2}$, then the
convergence $\mathcal{L}\mbox{-}\lim_{n \rightarrow
\infty}(\Psi^n,Z^n) = (\Psi,B^H)$ implies the convergence
\[
    \mathcal{L}\mbox{-}\lim_{n \rightarrow \infty}
    \left( \Psi^n,Z^n,\int \Psi^n dZ^n \right) =
    \left( \Psi,B^H, \int \Psi d B^H \right).
\]
\end{thm}

As an immediate corollary to Theorem \ref{thm3} we see that the
fluctuations of the price process \ref{YepsN} around its first
order approximation converge in distribution to a stochastic
integral with respect to fractional Brownian motion.

\begin{cor} Let $\Psi$ be a continuous semi martingale with Doob-Meyer
decomposition $\Psi= M+A$. If $\mathbb{E} \{[M,M]_{T}\}<\infty$,
$\mathbb{E} \{|A|_{T}\}<\infty$ and $\mu \neq 0$, then there
exists $\sigma > 0$ such that
\begin{equation} %\label{limit1}
    \mathcal{L}\mbox{-}\lim_{T \rightarrow \infty} \mathcal{L}\mbox{-}\lim_{N
    \rightarrow \infty} T^{1-H}\frac{\sqrt{N}}{\sqrt{L(T)}}
    \left\{ S^{N,T}_t - \mu \int_0^t \Psi_s ds \right\}_{0 \leq t \leq 1}
    =  \left\{ \sigma \int_0^t \Psi_s dB^H_s \right\}_{0 \leq t \leq 1}.
\end{equation}
%For the special case $\Psi \equiv 1$ the fluctuations can be
%described by a fractional Brownian motion.
\end{cor}

The increments of a \fbm with Hurst coefficient $H\in(\half,1]$
are positively correlated. The correlation increases in $H$. Thus,
the limit theorem reveals that, in isolation, investor inertia may
lead to long range dependence in asset returns. Indeed, a greater
degree of inactivity, represented by a smaller tail index
$\alpha$, leads to a larger $H$, and so greater positive
correlation between returns. Since \fbm is not a semimartingale,
it may also lead to arbitrage opportunities for other traders
whose impact has not been considered in the model so far. Explicit
arbitrage strategies for various models were constructed in, e.g.
\cite{erhan}.

\begin{rem}
In a model without inertia where all the sojourn time
distributions are thin-tailed, the logarithmic stock price
fluctuations can be approximated in law by a process of the form
\begin{equation} \label{day-traders}
    \left\{ \int_0^t \Psi_s\, dW_s \right\}_{0 \leq t \leq 1}
\end{equation}
where $W$ is a standard Brownian motion. Thus, when all traders'
mood processes are standard Markov processes and $\Psi$ is
constant, we recover in the limit the standard
Black-Scholes-Samuelson geometric Brownian motion model.
\end{rem}

The approach of studying queuing systems through their limiting
behaviour has a long history in many applications, see
\cite{Whitt}, for example. This analysis of investor inertia built
upon the works of Taqqu {\em et al.} \cite{taqqu} on internet
traffic. However, even the simple model we have discussed so far
shows how economic applications lead to new mathematical
challenges: in the teletraffic application, it is sufficient to
consider a binary (on/off) state space, but when agents buy, sell
or do nothing, there must be at least three states. This requires
different techniques from the binary case. Our functional central
limit theorems for stationary semi-Markov processes may also serve
as a mathematical basis for proving heavy-traffic limits in the
multilevel network models studied in, e.g. \cite{Duffield-Whitt}
and \cite{D-W}.

\subsection{A Limit Theorem with Feedback Effects} \label{Sec-Feedback}

The model in the previous section assumes that investors' actions
affect the price, but prices did not affect the agents' demands.
This assumption might be justified for Internet or new economy
stocks where no accurate information about the actual underlying
fundamental value is available. In such a situation, price is not
always a good indicator of value and is often ignored by
uninformed small investors. In general, however, it is certainly
desirable to allow for feedback effects from current prices into
the agents' order rates. In this section we extend our previous
model to allow for feedback effects from prices into the agents'
order rates. At the same time we provide a unified mathematical
framework for analyzing microstructure models with asynchronous
order arrivals. Our approach is based on methods and techniques
from state dependent Markovian service networks. Mathematically,
it extends earlier results in \cite{Anisimov} beyond semi-Markov
models with thin-tailed sojourn time distributions.

\subsubsection{The dynamics of logarithmic asset prices}\label{sec:our-feedback-model}
Let us now be more precise about the probabilistic structure our
model. We assume that the agents' orders arrive with an order rate
that depends on the price and the investor sentiment. Each order
is good for one unit of the stock. Specifically, we associate to
each agent $a \in \ma$ two independent standard Poisson processes
$\left\{ \Pi_+^a(t) \right\}_{t \geq 0}$ and $\left\{ \Pi_-^a(t)
\right\}_{t \geq 0}$, a stationary semi-Markov process $x^a$ on
$E$ satisfying Assumption~\ref{Assumption-tails}, and bounded
Lipschitz continuous \textsl{rate functions} $\lambda_\pm : E
\times \re \rightarrow \re^+$. The rate functions along with the
Poisson processes $\Pi^a_\pm$ specify the arrivals times of buying
and selling orders. The agent's holdings at time $t \geq 0$ are
given by

\begin{equation}\label{eq:an-agent}
    \Pi^a_+\left(\int_0^t\lambda_+\left( x^a_u, S^N_u \right)\,du \right)
    -
    \Pi^a_-\left(\int_0^t\lambda_-\left(x^a_u,S^N_u\right)\,du\right)
\end{equation}
where $\{S^N_t\}_{t \geq 0}$ denotes the logarithmic asset price
process. As before, a buying (selling) order increases (decreases)
the logarithmic price by $1/N$. Assuming for simplicity that
$S^N_0 = 0$, we thus obtain
\begin{equation} \label{eq:mim}
    S^N_t = \frac{1}{N} \sum_{a \in \ma}
    \Pi^a_+\left(\int_0^t\lambda_+\left( x^a_u, S^N_u \right)\,du \right)
    - \frac{1}{N} \sum_{a \in \ma}
    \Pi^a_-\left(\int_0^t\lambda_-\left(x^a_u,S^N_u\right)\,du\right).
\end{equation}

\begin{rem}
\begin{itemize}
\item[(i)] In the model studied in the previous section, the
agents continuously accumulated the stock at rates specified by
semi-Markov processes. Our current models assume that stocks are
purchased at random points in times. The arrival times of buying
and selling times follow exponential distributions conditional on
random arrival rates that depend on current prices and exogenous
semi-Markov processes.

\item[(i)] As before, we think of $x^a$ as being the investor's
``mood'' (for trading) process. Loosely speaking,
$\lambda_+(x^a_t,s) - \lambda_-(x^a_t,s)$ may be viewed as the
agent's excess demand at time $t$ at a logarithmic price level
$s$, given his trading mood $x^a_t$.

\item[(ii)] To develop a model of interaction, in which the
participants are inert, out of (\ref{eq:an-agent}), it is natural
to assume that $\lambda_\pm(0,s) \equiv 0$ and that the buying and
selling rates $\lambda_+(x,\cdot)$ and $\lambda_-(x,\cdot)$ are
increasing, rep. decreasing, in the second variable meaning that
meaning high (low) prices temper buying (selling) rates.
\end{itemize}
\end{rem}

The sum of independent Poisson processes is a Poisson process with
intensity given by the sum of the intensities. As a result, the
logarithmic price process satisfies the equality
\begin{equation}\label{eq:market-imbalance}
    S^N_t = \frac{1}{N} \Pi_+\left ( \sum_{a = 1}^N \int_0^t
    \lambda_+ \left(x^a_u, S^N_u \right) \,du
    \right) - \frac{1}{N} \Pi_-\left ( \sum_{a = 1}^N \int_0^t
    \lambda_-\left(x^a_u,S^N_u \right)\,du \right)
\end{equation}
in \textsl{distribution} where $\Pi_+$ and $\Pi_-$ are independent
standard Poisson processes. Since our focus will be on a limit
result for the \textsl{distribution} of the price process as the
number of agents grows to infinity, we may with no loss of
generality assume that the logarithmic price process is
\textsl{defined} by (\ref{eq:market-imbalance}) rather than
(\ref{eq:mim}).

\begin{Assumption}\label{assumplambda}
\begin{enumerate}

\item The rate functions $\lambda_{\pm}$ are  uniformly bounded.

\item  For each $x \in E$, the rate functions $\lambda_{\pm}(x,
\cdot)$ are continuously differentiable with first derivative
bounded in absolute value by some constant $L$.
\end{enumerate}
\end{Assumption}

Our convergence results will be based on the following strong
approximation result which allows for a pathwise approximation of
a Poisson process by a standard Brownian motion living on the same
probability space.

\begin{lemma} (\cite{kurtz}) \label{strong-approximation}
    A standard Poisson process $\left\{\Pi(t)\right\}_{t \geq 0}$ can be
    realized on the same probability space as a standard Brownian motion
    $\{B(t)\}_{t \geq 0}$ in such a way that the almost surely finite random variable
\[
    \sup_{t \geq 0} \frac{|\Pi(t) - t - B(t)|}{\log (2 \vee t)}
\]
    has a finite moment generating function in the neighborhood of
    the origin and in particular finite mean.
\end{lemma}

In view of Assumption \ref{assumplambda} (i), the strong
approximation result yields the following alternative
representation of the logarithmic asset price process:
\begin{equation} \label{representation-S}
\begin{split}
    S^N_t = & ~ \frac{1}{N} \left\{ \sum_{a = 1}^N \int_0^t
    \lambda \left(x^a_u, S^N_u \right)du
    +
    B_+ \left( \sum_{a = 1}^N \int_0^t \lambda_+
    \left(x^a_u, S^N_u \right)du \right) \right. \\
    & ~ \left. - B_- \left( \sum_{a = 1}^N \int_0^t \lambda_-
    \left(x^a_u, S^N_u \right)du \right) \right\} +
    {\cal O}\left( \frac{\log N}{N} \right),
\end{split}
\end{equation}
where $\lambda(x^a_u,\cdot)$ denotes the excess order rate of the
agent $a \in \ma$, given his mood for trading $x^a_u$ and ${\cal
O}\left( \log N / N \right)$ holds uniformly over compact time
intervals. Using this representation of the logarithmic price
process our goal is to prove approximation results for the process
$\{S^N_t\}_{t \geq 0}$. In a first step we show that it can almost
surely be approximated by the trajectory of an ordinary
differential equation (``fluid limit''). In subsequent step, we
apply a result from \cite{kn:bhs} to show that, after suitable
scaling, the fluctuations around this first order approximation
can be described in terms of a fractional process $\{Z_t\}_{t \geq
0}$ of the form
\[
    d Z_t = \mu_t Z_t dt + \sigma_t d B^H_t.
\]
In a benchmark model without feedback, where the order rates do
not depend on current prices, the process $\{Z_t\}_{t \geq 0}$
reduces to a fractional Brownian motion. That is, we recover the
type of results of Section~\ref{sec:investor-inertia} with the
alternative model presented in this section.

%%%%%%%%%%%%%%%%%%%%%%%%%%%%%%%%%%%%%%%%%%%%%%%%%%%%%%%%%%%%%%%%%%%%%%%%%%

\subsubsection{First order approximation}

In order to prove our first convergence result, it is convenient
to denote by
\begin{equation}\label{eq:defn-lambda}
    \lambda(x,s) \triangleq  \lambda_+(x,s) - \lambda_{-}(x,s)
\end{equation}
the accumulated net order rate at a given logarithmic price level
$s \in \re$ and trading mood $x \in E$ and by
\[
    \bar{\lambda}(s) \triangleq  \bar{\lambda}_+(s) - \bar{\lambda}_{-}(s)
\]
the expected excess order flow where
\[
    \bar{\lambda}_{\pm}(s) \triangleq  \int_E \lambda_{\pm}(x,s)
    \nu
    (dx),
\]
and $\nu$ is the stationary distribution of the semi-Markov
process $x_t$. We are first going to show that in a financial
market with many agents the dynamics of the logarithmic price
process can be approximated by the solution $\{s_t\}_{t \geq 0}$
to the ODE
\begin{equation}\label{eq:defn-of-s}
    \frac{d}{dt}\,s_t = \bar{\lambda} (s_t),
\end{equation}
with initial condition $s_0=0$. To this end, we need to prove that
the average excess order rate converges almost surely to the
expected excess order flow uniformly on compact time intervals.

\begin{lemma}\label{lem:uniform-LLN}
Uniformly on compact time intervals
\begin{equation}
    \lim_{N \rightarrow \infty} \frac{1}{N} \sum_{a=1}^N \int_0^t
    \lambda_{\pm} (x^a_u,s_u) du = \int_0^t \bar{\lambda}_\pm (s_u) du \qquad
    \Prob^*\mbox{-a.s.}
\end{equation}
\end{lemma}
\begin{Proof}
The stationary semi Markov processes $x^a$ are independent, and so
the random variables $\int_{0}^{t}\lambda(x_{u}^{a},s_u)du$
$(a=1,2,...)$ are also independent. Thus, the law of large numbers
for independent random variables along with Fubini's theorem (to
exchange the sum and the integral) and bounded convergence theorem
(to exchange the limit and the integral) yields convergence for
each $t$. In order to prove that the convergence holds uniformly
over compact time intervals we will use uniform law of large
numbers of \cite{prucha}. Denoting $\textnormal{D}_E[0,t]$ the
class of all c\'{a}dl\'{a}g functions $y:[0,t] \rightarrow E$ we
need to show that the maps $q_\pm: \textnormal{D}_E[0,t] \times
[0,t] \rightarrow \re$ defined by
\[
    q_\pm(y,t) \triangleq \int_0^t \lambda_\pm(y(u),s_u)du
\]
are continuous. Since the rate functions are bounded, it is enough
to show that the map $y \mapsto \int_0^\cdot \lambda_\pm(y(u),s_u)
du$ is continuous uniformly over compact time intervals.

To this end, we denote by $d$ the metric defined in (3.5.2) in
\cite{ethier} which induces the Skorohod topology in $D_E[0,t]$
and recall that $\lim_{n \rightarrow \infty} d(y_n,y) = 0$ if and
only if
\begin{equation} \label{convergence-D}
    \lim_{n \rightarrow \infty} \sup_{0 \leq s \leq
    t}|y_n \circ \tau_n(s) - y(s)| = 0
\end{equation}
for a suitable sequence of strictly increasing time-shifts
$\tau_n$; see \cite{ethier} page 117 for details. Let $\{y_n\}$
denote a sequence in $D_E[0,t]$ that converges to $y$ and put
\[
    \lambda_\pm^n(u) \triangleq \lambda_\pm(y_n(u),s_u).
\]
In view of the transformation formula for Lebesgue integrals and
because $\tau(0) = 0$ and $\tau_n^{-1}(t) \leq t$ we obtain
\begin{eqnarray*}
    \int_0^t \left[ \lambda^n_\pm(u) - \lambda_\pm(u) \right] du
    & = & \int_0^{\tau_n^{-1}(t)} \left[
    \lambda^n_\pm \circ \tau_n(u) \tau'_n(u) - \lambda_\pm(u) \right]du -
    \int_{\tau_n^{-1}(t)}^t \lambda_\pm(u) du \\
    & = & \int_0^{\tau_n^{-1}(t)} \left[
    \lambda^n_\pm \circ \tau_n(u) - \lambda_\pm(u) \right]du \\ & & +
    \int_0^{\tau_n^{-1}(t)} \lambda^n_\pm \circ \tau_n(u)[\tau'_n(u) - 1] du
    - \int_{\tau_n^{-1}(t)}^t \lambda_\pm(u) du.
\end{eqnarray*}
    By (3.5.5)-(3.5.7) in \cite{ethier}
\[
    \lim_{n \rightarrow \infty} \sup_{0 \leq u \leq
    t} |\tau_{n}'(u) - 1| = 0 \quad \mbox{and} \quad
    \lim_{n \rightarrow \infty} \sup_{0 \leq u \leq
    t} |\tau^{-1}_{n}(u) - u| = 0
\]
so that the last two terms on the right hand side of the
inequality above vanish uniformly on compact time intervals. As
far as the first term is concerned, observe that boundedness of
the rate function's derivative with respect to the second argument
yields
\[
    \left| \lambda_\pm(y_n \circ \tau_n (u), s_{\tau_n(u)}) -
    \lambda_\pm(y(u),s_u) \right| \leq L \left|y_n \circ \tau_n (u) -
    y(u) \right| + L|s \circ \tau_n(u) - s(u) |.
\]
    As a continuous function $s$ is uniformly continuous over compact time
    intervals. This, along with (\ref{convergence-D}) yields
\[
    \lim_{n \rightarrow \infty} \sup_{0 \leq u \leq t}
    \left| \lambda_\pm(y_n \circ \tau_n (u), s_{\tau_n(u)}) - \lambda_\pm(y(u),s_u)
    \right| = 0
\]
    so that the maps $q_\pm$ are indeed continuous. Thus, the uniform
    law of large numbers yields
\[
    \lim_{N \rightarrow \infty}
    \frac{1}{N} \sum_{a=1}^N q_\pm (x^a_{u},u) = \lim_{N \rightarrow \infty}
    \frac{1}{N} \sum_{a=1}^N \int_0^u \lambda_\pm(x^a_v,s_v)dv =
    \lambda_\pm(\mu,s_u) %\quad \Prob\mbox{-a.s.}
\]
    almost surely uniformly on compact time intervals.
\end{Proof}

We are now ready to state and prove our functional law of large
numbers.

\begin{thm}\label{lem:fluidlimit}
    As $N \rightarrow \infty$, the sequence of stochastic processes
    $\{S^N_t\}_{t \geq 0}$ $(N \in \n)$ converges almost surely to the
    deterministic process $\{s_t\}_{t \geq 0}$:
\[
    \lim_{N \rightarrow \infty} S^N_t = s_t \qquad
    \Prob^*\mbox{-a.s.}
\]
    where the convergence is uniform over compact time intervals.
\end{thm}
\begin{Proof}
In view of the strong approximation result formulated in Lemma
\ref{strong-approximation} and because the rate functions are
uniformly bounded,
\begin{eqnarray*}
    \left| \Pi_\pm \left ( \sum_{a = 1}^N \int_0^t \lambda_\pm
    \left(x^a_u, S^N_u\right)du \right) - \sum_{a = 1}^N \int_0^t
    \lambda_\pm \left(x^a_u, S^N_u\right)du - B_\pm \left( \sum_{a = 1}^N
    \int_0^t \lambda_\pm \left(x^a_u, S^N_u\right)du \right) \right|
\end{eqnarray*}
is of the order $O \left( \log N \right)$ almost surely where
$B_{\pm}$ are the Brownian motions used in
(\ref{representation-S}). Since the rate functions are uniformly
bounded, the law of iterated logarithm for
Brownian motion yields %that
%\[
%    \limsup_{N \rightarrow \infty}
%    \frac{B_\pm \left( \sum_{a = 1}^N \int_0^u \lambda_\pm
%    \left(x^a_v, S^N_{v} \right)dv \right)}
%    {\left(2 N C u \log \log (NCu) \right)^{1/2}} \leq 1 \qquad
%    \Prob^*\mbox{-a.s.}
%\]
%for any $u \geq 0$. Therefore
\[
    \lim_{N \rightarrow\infty} \sup_{u \leq t} \frac{1}{N}
    B_\pm \left( \sum_{a = 1}^N \int_0^u \lambda_\pm
    \left(x^a_v, S^{N}_v \right)dv \right) = 0 \qquad \Prob^*\mbox{-a.s.}
\]
It follows from this and Lemma \ref{lem:uniform-LLN} above, that
the quantities
\[
    B^N_t \triangleq  \frac{1}{N} \left|  B_+ \left( \sum_{a = 1}^N \int_0^t
    \lambda_+ \left(x^a_u, S^N_u \right)du \right) -
    B_- \left( \sum_{a = 1}^N \int_0^t \lambda_-
    \left(x^a_u, S^{N}_u \right)du \right) \right|
\]
and
\[
    \Lambda^N_t \triangleq \left| \frac{1}{N} \sum_{a  = 1}^N \int_0^t
    \left\{ \lambda \left( x^a_u, s_u \right) -
    \bar{\lambda} \left(s_u \right)
    \right\} du \right|
\]
converge to zero uniformly over compact time intervals as $N
\rightarrow \infty$.

Let us now fix $\epsilon > 0$. Due to Lemma
\ref{strong-approximation} there exists $N^* \in \n$ such that for
all $N \geq N^*$ and uniformly on compact time sets, for $l \leq
t$ we can write
\begin{eqnarray*}
    \left| S^N_l- s_l \right| & \leq &
    \left| \frac{1}{N} \sum_{a  = 1}^N \int_0^l
    \lambda \left( x^a_u,  S^N_u \right) du -
    \int_0^l \bar{\lambda}(s_u) du \right| + B^N_l + \epsilon \\
    & \leq & \left| \frac{1}{N} \sum_{a  = 1}^N \int_0^l
    \left\{
    \lambda \left( x^a_u, S^N_u \right) -
    \lambda \left( x^a_u, s_u \right)
    \right\} du  \right| + \Lambda^N_l + B^N_l + \epsilon \qquad
    \Prob^*\mbox{-a.s.}
\end{eqnarray*}
Lipschitz continuity of the rate functions yields
\begin{eqnarray*}
    \left| S^N_l - s_l \right|
    & \leq & L \int_0^l \sup_{0 \leq r \leq u} \left |
    S^N_r - s_r \right| du + \Lambda^N_l + B^N_l +
    \epsilon \\
     & \leq & L \int_0^t \sup_{0 \leq r \leq u} \left |
    S^N_r - s_r \right| du
    + \sup_{0 \leq r \leq t} \Lambda^N_r +
    \sup_{0 \leq r \leq t} B^N_r + \epsilon \qquad \Prob^*\mbox{-a.s.}
\end{eqnarray*}
for some $L>0$ and so
\begin{equation} \label{Gronwall-S}
    \sup_{0 \leq r \leq t} \left| S^N_r - s_r \right|
    \leq L \int_0^t \sup_{0 \leq r \leq u} \left |
    S^N_r - s_r \right| du + \sup_{0 \leq r \leq t} \Lambda^N_r +
    \sup_{0 \leq r \leq t} B^N_r + \epsilon \qquad \Prob^*\mbox{-a.s.}
\end{equation}
Now, an application of Gronwall's lemma yields
\[
    \sup_{0 \leq r \leq t} \left| S^N_r - s_r \right|
    \leq \left( \sup_{0 \leq r \leq t} \Lambda^N_r +
    \sup_{0 \leq r \leq t} B^N_r + \epsilon \right) e^{L t}
    \qquad \Prob^*\mbox{-a.s.}
\]
for all $N \geq N^*$. This proves our assertion.
\end{Proof}

%%%%%%%%%%%%%%%%%%%%%%%%%%%%%%%%%%%%%%%%%%%%%%%%%%%%%%%%%%%%%%%%%%%%%%%%%%%%%%%%%%%%%

\subsubsection{Second order approximation}

In this section we analyze the fluctuations of the logarithmic
price process around its first order approximation. We are
interested in the distribution of asset prices around their first
order approximation as $N \rightarrow \infty$. In view of the
representation (\ref{representation-S}) and by self-similarity of
Brownian motion we may thus assume that $\{S^N_t\}_{t \geq 0}$ is
defined by the integral equation:
\begin{equation} \label{price2}
\begin{split}
    S^N_t =~ & \frac{1}{N} \sum_{a = 1}^N \int_0^t
    \lambda \left(x^a_u, S^N_u\right)du
    + \frac{1}{\sqrt{N}} B_+ \left( \frac{1}{N}
    \sum_{a = 1}^N \int_0^t \lambda_+
    \left(x^a_u, S^N_u \right)du \right) \\
    & ~ - \frac{1}{\sqrt{N}} B_- \left( \frac{1}{N}
    \sum_{a = 1}^N \int_0^t \lambda_-
    \left(x^a_u, S^N_u\right)du \right) + O \left(\frac{\log N}{N} \right).
\end{split}
\end{equation}
As we shall see, the fluctuations around the first order
approximation are driven by two Gaussian processes. The first,
\begin{equation} \label{process-X}
    X_t \triangleq B_+ \left( \int_0^t \bar{\lambda}_+ (s_u)du \right) -
    B_- \left( \int_0^t \bar{\lambda}_- (s_u)du \right),
\end{equation}
captures the randomness in the agents' trading times. The second,
$\{Y_t\}_{t \geq 0}$, is defined in terms of the integral of a
non-stationary Gaussian process whose covariance function depends
on the first order approximation. It captures the second source
randomness generated by the agents' trading activity.
Specifically,
\begin{equation}\label{eq:YT}
    Y_t \triangleq \int_0^t y_s ds,
\end{equation}
where $\{y_t\}_{t \geq 0}$ denotes the centered Gaussian process
whose covariance function $\gamma$ is given by the covariance
function of the stochastic process $\{\lambda(x_t,s_t)\}_{t \geq
0}$, i.e.,
\begin{equation} \label{gamma}
    \gamma(t,u) \triangleq \E[\lambda(x_t,s_t) \lambda(x_u,s_u)] -
    \bar{\lambda}(s_t) \bar{\lambda}(s_u).
\end{equation}
It turns out that the fluctuations can be approximated in
distribution by the process $\{Z_t\}_{t \geq 0}$ which satisfies
the integral equation %defined by
%\begin{equation}
%    Z_{t} = \int_{0}^{t}\exp\left(\int_{s}^{t}\bar{\lambda}'(s_{u})du\right)
%    \, (d X_{s}+ d Y_s)
%\end{equation}
%if the number os agents grows to infinity. Note that $\{Z_t\}_{t
%\geq 0}$ satisfies the integral equation
\begin{equation} \label{Z}
    Z_t = \int_0^t \bar{\lambda}'(s_u) Z_u du + Y_t + X_t.
\end{equation}

Our goal is to establish the following second order approximation
for the asset price process in an economy with many market
participants. %Without loss of generality we will assume that the
%price process evolves on the unit interval $[0,1]$.

\begin{thm} \label{thm-limit-N}
    The fluctuations of the market imbalance
    $\{S^N_t\}_{0 \leq t \leq 1}$ around its first order approximation can be
    described by the process $\{Z_t\}_{0 \leq t \leq 1}$ defined in (\ref{Z}).
    More precisely,
\[
    \mathcal{L}\mbox{-}\lim_{N \rightarrow \infty}
    \sqrt{N}\left\{ S^N_t - s_t \right\}_{0 \leq t \leq 1}
    = \left\{ Z_t \right\}_{0 \leq t \leq 1}.
\]
\end{thm}

\mbox{ }

The proof of Theorem \ref{thm-limit-N} requires some preparation.
For notational convenience we introduce stochastic processes
$Q^N=\{Q^N_t\}_{0 \leq t \leq 1}$, $Y^N = \{Y^N_t\}_{0 \leq t \leq
1}$ and $X^N = \{X^N_t\}_{0 \leq t \leq 1}$ by, respectively,
\begin{equation}\label{eq:defn-YN}
    Q^N_t \triangleq \sqrt{N} (S^N_t - s_t) \qquad \mbox{and}
    \qquad
    Y^N_t \triangleq  \sum_{a=1}^N \int_0^t \frac{ \lambda(x^a_u,s_u) -
    \bar{\lambda}(s_u)}{\sqrt{N}} \ du,
\end{equation}
and
\begin{equation} \label{ZN}
    X^N_t  \triangleq  B_+
    \left( \frac{1}{N} \sum_{a = 1}^N \int_0^t \lambda_+ \left(x^a_u, S^N_u \right)du
    \right) - B_-
    \left( \frac{1}{N} \sum_{a = 1}^N \int_0^t \lambda_- \left(x^a_u, S^N_u\right)du
    \right).
\end{equation}
We first prove convergence in distribution of the sequence
$\{(X^N,Y^N)\}_{N \in \n}$ to $(X,Y)$.

\begin{prop} \label{lemma-joint-convergence}
The sequence $\{(X^N, Y^N)\}_{N \in \n}$ converges in distribution
to the process $(X, Y)$ defined by (\ref{process-X}) and
(\ref{eq:YT}).
\end{prop}
\begin{Proof} For any $\alpha \in (0,\frac{1}{2})$ and
$T>0$, there exist integrable and hence almost surely finite
random variables $M_\pm$ such that for all $t_1,t_2 \leq T$ we
have
\[
    |B_\pm(t_1) - B_\pm(t_2)| \leq M_\pm |t_1 - t_2|^\alpha \quad
    \Prob^*\mbox{-a.s.},
\]
see, for instance, Remark 2.12 in \cite{kn:karat}. Thus, the first
order approximation shows that the sequence of processes
$\{X^N\}_{N \in \n}$ converges almost surely to $X$ on any compact
time interval. Since the processes
\[ \int_0^t \frac{ \lambda(x^a_u,s_u) -
    \bar{\lambda}(s_u)}{\sqrt{N}} \ du
\]
have Lipschitz continuous sample paths and the semi-Markov
processes are independent, the central limit theorem for Lipschitz
processes (\cite{Whitt}, Corollary 7.2.1) shows that $\{Y^N\}_{N
\in \n}$ converges in distribution to the Gaussian process $Y$. As
a result, both sequences $\{X^N\}_{N \in \n}$ and $\{Y^N\}_{N \in
\n}$ are tight. Since $\{X^N\}_{n \in \n}$ is also C-tight, the
sequence $\{(X^N, Y^N)\}_{N \in \n}$ is tight. It is therefore
enough to prove weak convergence of the finite dimensional
distributions of the process $(X^N, Y^N)$ to the finite
dimensional distributions of $(X, Y)$.

In order to establish weak convergence of the one-dimensional
distributions we fix a Lipschitz continuous functions with compact
support $F: \re^2 \rightarrow \re$. We may with no loss of
generality assume that both the Lipschitz constant and the
diameter of the support of $F$ equal one. In this case
\[
    \left|
    \int F \left(X^N_t, Y^N_t \right) d
    \Prob^* - \int F \left(X_t, Y^N_t \right) d
    \Prob^* \right| \leq \int \min \{|X^N_t - X_t|,1\} d
    \Prob^*.
\]
In view of the convergence properties of the sequence $\{X^N\}_{N
\in \n}$, there exists, for any $\epsilon > 0$, a constant $N^*
\in \n$ such that
\[
    \sup_{0 \leq t \leq 1} \int \min \{|X^N_t - X_t|,1\} d
    \Prob^* \leq \epsilon \quad \mbox{for all} \quad N \geq N^*.
\]
This yields
\[
    \lim_{N \rightarrow \infty} \left|
    \int F \left(X^N_t, Y^N_t \right) d
    \Prob^* - \int F \left(X_t, Y^N_t \right) d
    \Prob^* \right| = 0.
\]
Since the random variables $X_t$ and $Y^N_t$ are independent, we
also have that
\[
    \lim_{N \rightarrow \infty }\int F \left(X_t,Y^N_t \right)
    d \Prob^* = \int F\left(X_t, Y_t \right) d \Prob^*.
\]
This proves vague convergence\footnote{A sequence of probability
measure $\{\mu_n\}$ converges to a measure $\mu$ in the vague
topology if $\lim_{n \rightarrow \infty} \int f d \mu_n = \int f
d\mu$ for all continuous functions $f$ with bounded support. The
vague limit $\mu$ is not necessarily a probability measure.
However, if there is an {\em a priori} reason that $\mu$
\textsl{is} a probability measure, then weak convergence of
$\{\mu_n\}$ to $\mu$ can be established by analyzing integrals of
continuous and hence Lipschitz continuous functions with bounded
support. See e.g. \cite{bauer} and \cite{Billingsley2}.} of the
one-dimensional marginal distributions of $(X^N,Y^N)$ to the
one-dimensional distributions of $(X,Y)$ and hence weak
convergence. Weak Convergence of the finite dimensional
distributions follows from similar considerations.
\end{Proof}

The following ``compact containment condition'' is key to the
second order approximation.
% to the
%representation (\ref{representation-integral}).

\begin{lemma} \label{lemma-QN-tight}
\begin{rmenumerate}
\item The sequence of stochastic processes $\{Q^N\}_{N \in \n}$ is
bounded in probability. That is, for any $\epsilon > 0$, there
exists $N^* \in \n$ and $K < \infty$ such that
\begin{equation} \label{tightness-QN}
    \Prob^*\left[ \sup_{0 \leq t \leq 1} |Q^N_t| > K \right] <
    \epsilon \quad \mbox{for all} \quad N \geq N^*.
\end{equation}
\item If $f^N = \{f^N_t\}_{t \geq 0}$ be a sequence of
non-negative random processes such that
\begin{equation}
    \lim_{N \rightarrow \infty} \int_0^1 f^N_{u} du = 0 \quad
    \mbox{in probability,}
\label{fcond}\end{equation}
    then, for all $\delta > 0$,
\[
    \lim_{N \rightarrow \infty}
    \Prob^* \left[ \sup_{0 \leq t \leq 1} \left| \int_0^t Q^N_u f^N_{u} du \right|
  > \delta \right] = 0.
\]
\end{rmenumerate}
\end{lemma}
\begin{Proof}
\begin{rmenumerate}
\item The strong approximation for Brownian motion yields the
representation
\begin{equation} \label{representation-QN}
    Q^N_t = \frac{\int_0^t \sum_{a=1}^N
    \left\{ \lambda\left(x^a_u, S^N_u\right) - \lambda(x^a_u,s_u)\right\} du}
    {\sqrt{N}} +  Y^N_t + X^N_t + O \left(
    \frac{\log N}{\sqrt{N}}\right).
\end{equation}
By Proposition \ref{lemma-joint-convergence} the sequence
$\{(X^N,Y^N)\}_{n \in \n}$ is tight, and hence it is bounded in
probability (see e.g. \cite{D-W}). As a result, Lipschitz
continuity of the rate functions yields
\[
    \sup_{0 \leq t \leq 1} |Q^N_t| \leq L \int_0^T \sup_{0 \leq t \leq
    u} |Q^N_u| du + \sup_{0 \leq t \leq 1} \left| Y^N_t \right| +
    \sup_{0 \leq t \leq 1} |X^N_t| + {\cal O}\left( \frac{\log N}{\sqrt{N}} \right).
\]
for some $L>0$. Hence, by Gronwall's inequality,
\[
    \sup_{0 \leq t \leq 1} |Q^N_t| \leq e^{LT} \left[
    \sup_{0 \leq t \leq 1} \left| Y^N_t \right| +
    \sup_{0 \leq t \leq 1} |X^N_t| + {\cal O}\left( \frac{\log N}{\sqrt{N}}
    \right) \right] \quad \Prob^*\mbox{-a.s.}
\]
This proves (i).

\item Let us fix $\epsilon>0$. There exists a constant $N^{*}$
such that when $N \geq N^{*}$ there exist sets $\Omega_N$ and
$A_N$ such that
\[
    \int_0^1 f^N_{u} du < \frac{\epsilon}{2} \quad \mbox{on $\Omega_N$}
    \quad \mbox{and such that} \quad \Prob^*[\Omega_N] \geq 1 -
    \frac{\epsilon}{2}.
\]
and
\[
    \sup_{0 \leq t \leq 1} |Q^N_t| < K \quad \mbox{on $A_N$}
    \quad \mbox{and such that} \quad \Prob^*[A_N] \geq 1 -
    \frac{\epsilon}{2}.
\]
Hence
\[
    \sup_{0 \leq t \leq 1} \left| \int_0^t Q^N_u f^N_u du \right| \leq \sup_{0 \leq t \leq
    1} |Q^N_t| \int_0^1 f^N_{u} du < K \epsilon \quad \mbox{on $A_N \cap \Omega_N$.}
\]
\end{rmenumerate}
\end{Proof}

\mbox{ }

\noindent \textsl{Proof of Theorem \ref{thm-limit-N}:} Let us
first define a sequence of stochastic processes $\widetilde{Q}^N =
\{\widetilde{Q}^N_t\}_{0 \leq t \leq 1}$ by
\[
    \widetilde{Q}^N_t \triangleq \int_0^t \bar{\lambda}'(s_u)
    \widetilde{Q}^N_u du + Y^N_t + X^N_t.
\]
By the continuous mapping theorem and Lemma
\ref{lemma-joint-convergence} the sequence $\{\widetilde{Q}^N\}_{N
\in \n}$ converges in distribution to the process $Z$ defined in
(\ref{Z}). It is now enough to show that
\begin{equation} \label{probabilistic-limit}
    \lim_{N \rightarrow \infty} \sup_{0 \leq t \leq 1} |Q^N_t -
    \widetilde{Q}^N_t | = 0 \quad \mbox{in probability.}
\end{equation}
To this end, let $E^{N}_{t} \triangleq Q^{N}_{t}-\tilde{Q}^{N}_t$.
>From the definition of $\tilde{Q}^{N}_t$ and the representation
(\ref{representation-QN}) of $Q^N_t$ we obtain

\begin{eqnarray*}
    E^{N}_{t} & = & \int_{0}^{t} \bar{\lambda}'(s_u) E^{N}_u du+
    \frac{1}{\sqrt{N}} \int_0^t \sum_{a=1}^N
    \bigg \{ \lambda\left(x^a_u, S^N_u \right)
    - \lambda(x^a_u,s_u) \bigg\} du
    - \int_{0}^{t}\bar{\lambda}'(s_u)Q^{N}_{u}du \\ %+ O \left(\frac{\log N}{\sqrt N}
    % \right)
    & = & \int_{0}^{t} \bar{\lambda}'(s_u) E^{N}_u du + \int_0^{t}
    \left(\frac{1}{N}\sum_{a=1}^{N}
    \lambda'(x^{a}_u,s_u)-\bar{\lambda}'(s_u)\right)Q^{N}_u\, du
    \\
    & & + \int_{0}^{t} \left(\frac{1}{N}\sum_{a=1}^{N}
    \lambda'(x^{a}_u,\xi^N_u) - \lambda'(x^{a}_u,s_u) \right)
    Q^{N}_{u}\,du. %+ O\left(\frac{\log N}{\sqrt N}\right)
\end{eqnarray*}
The second equality follows from the mean value theorem for
$\lambda(x^a_u,\cdot)$,
\[
    \lambda\left(x^a_u, S^N_u \right)-\lambda(x^a_u,s_u)
    = \frac{1}{\sqrt{N}}\lambda'(x^{a}_u, \xi^{N}_u) Q^{N}_u,
\]
where $\xi^{N}_u$ lies between $\frac{1}{N}S^{N}_u$ and $s_u$. We
put
\[
    f^{N,1}_u \triangleq \frac{1}{N}\sum_{a=1}^{N}  \lambda'(x^{a}_u,s_u)
    - \bar{\lambda}'(s_u) \quad \mbox{and} \quad
    f^{N,2}_u \triangleq \frac{1}{N}\sum_{a=1}^{N}
    \lambda'(x^{a}_u,\xi^N_u)-\lambda'(x^{a}_u,s_u)
\]
in order to obtain
\[
    \sup_{0 \leq s \leq t} |E^N_s| \leq L \int_0^t \sup_{0 \leq s \leq u} |E^N_s|\,du +
    \bigg| \sup_{0 \leq s \leq t} \int_0^s  |f^{N,1}_u | Q^N_u \, du \bigg| +
    \bigg| \sup_{0 \leq s \leq t} \int_0^s  |f^{N,2}_u| Q^N_u \, du\bigg|. %% +
    %% O \left(\frac{\log N}{\sqrt N} \right).
\]
The processes $|f^{N,1}|$ and $|f^{N,2}|$ satisfy the condition
(\ref{fcond}) of Lemma \ref{lemma-QN-tight} by the law of large
numbers. Thus, an application of Gronwall's lemma yields
(\ref{probabilistic-limit}). \hfill $\Box$

%%%%%%%%%%%%%%%%%%%%%%%%%%%%%%%%%%%%%%%%%%%%%%%%%%%%%%%%%%%%%%%%%%%

\subsubsection{Approximation by a fractional Ornstein-Uhlenbeck
process}

So far we have shown that the fluctuations of the logarithmic
price process around its first order approximation can be
described in terms of an Ornstein-Uhlenbeck process $Z$ driven by
two Gaussian processes $X$ and $Y$. In order to see more clearly
the effects of investor inertia on asset processes we need to
better understand the dynamics of $Y$. As before, this will be
achieved by a proper scaling of of the semi-Markov processes $x^a$
in time and the price process in space. Specifically, we introduce
a family of processes $S^{N,T}$ $(T \in \n)$ with initial value
$0$ by
\begin{eqnarray*}\label{eq:mimT}
%\begin{split}
    S^{N,T}_{t} = \frac{1}{NT} \left\{
    \Pi_+ \left(T \sum_{a=1}^N \int_0^{t} \lambda_+
    \left(x^a_{Tu}, S^{N,T}_u \right)\,du\right)
    -
    \Pi_-\left(T \sum_{a=1}^N \int_0^{t} \lambda_-
    \left(x^a_{Tu}, S^{N,T}_u \right)\,du
    \right) \right\}.% \\
\end{eqnarray*}
The strong approximation result for Poisson processes with respect
to Brownian motion allow us to represent the process
$\{S^{N,T}_{t}\}_{t \geq 0}$ as in (\ref{representation-S}) with
the semi-Markov processes $\{x^a_t\}_{t \geq 0}$ replaced by the
``speeded-up'' processes $\{x^a_{Tt}\}_{t \geq 0}$. Moreover, by
Lemma~\ref{lem:uniform-LLN}, the sequence of processes
\[
    \Lambda^{N,T}_t \triangleq \left| \frac{1}{N} \sum_{a  = 1}^N \int_0^t
    \left\{ \lambda \left( x^a_{Tu}, s_u \right) -
    \bar{\lambda} \left(s_u \right)
    \right\} du \right|
\]
converges to zero uniformly over compact time intervals as $N
\rightarrow \infty$. Following the same line of arguments as in
the proof of Proposition~\ref{lem:fluidlimit} it can then be shown
that for any $T > 0$
\begin{equation}
    \lim_{N \rightarrow \infty} S^{N,T}_{t} = s_t \qquad
    \Prob^*\mbox{-a.s.} \quad \mbox{}
\end{equation}
Here $\{s_t\}_{t \geq 0}$ denotes the deterministic process
defined by the ordinary differential equation (\ref{eq:defn-of-s})
with initial condition $s_0 = 0$ and the convergence holds
uniformly over compact time intervals. Thus, the first order
approximation is independent of $T$. By analogy to
(\ref{price2})-(\ref{Z}) introduce a Gaussian process $Y^T$ by
\begin{equation}
    Y^T_t \triangleq  \int_0^t y^T_s ds,
\end{equation}
where $\{y^T_t\}_{t \geq 0}$ denotes the centered Gaussian process
with covariance function
\[
%\left(\int_{0}^{t}\int_{0}^{s} \gamma^{T}(u,v)du dv\right)_{(u,v)
%\in \mathbb{R}^{2}_{+} }\quad \text{where} \quad
    \gamma^{T}(t,u) \triangleq
    \E[\lambda(x_{Tt},s_t) \lambda(x_{Tu},s_u)] -
    \bar{\lambda}(s_t) \bar{\lambda}(s_u).
\]
Following the same arguments in the proof of
Theorem~\ref{thm-limit-N}, we see that as the number of agents
tends to infinity the price fluctuations round the fluid limit can
be approximated in distribution by a process $\{Z^T_t\}_{t \geq
0}$ of the form
\[
    Z^T_t = \int_0^t \bar{\lambda}'(s_u) Z^T_u du + Y^T_t + \frac{1}{\sqrt{T}}X_t.
\]

\begin{prop}
    For any $T$, the fluctuations of the logarithmic price process
    $\{S^{N,T}_t\}_{0 \leq t \leq 1}$ around its first order approximation can be
    described by the process $\{Z^T_t\}_{0 \leq t \leq 1}$.
    More precisely,
\[
    \mathcal{L}\mbox{-}\lim_{N \rightarrow \infty}
    \sqrt{N}\left\{ S^{N,T}_{t} - s_t \right \}_{0 \leq t \leq 1}
    = \left\{ Z^T_t \right\}_{0 \leq t \leq 1}.
\]
\end{prop}

%
%Let us define
%\begin{equation}\label{eq:defn-Y-NT}
% Y^{N,T}_t \triangleq  \sum_{a=1}^N \int_0^t \frac{ \lambda(x^a_{Tu},s_u) -
%    \bar{\lambda}(s_u)}{\sqrt{N}} \ du
%\end{equation}
%As $N \rightarrow \infty$, $Y^{N,T}_t$ converges to a Gaussian
%process with covariance function

To take the $T$-limit, we need the following assumption on the
structure of the rate functions.

\begin{Assumption}\label{ass:separable}
The rate function $\lambda$ defined in (\ref{eq:defn-lambda}) can
be written as
\begin{equation}\label{eq:seperable}
\lambda(x,s)=f(x)g(s)+ h(s). %\quad f \neq 0,\,\, g \neq 0.
\end{equation}
Moreover, the function $f$ in (\ref{eq:seperable}) is one-to-one
and $\hat{\mu} \triangleq f(0) \neq \E^* f(x_0)$.
\end{Assumption}

\begin{example}
The previous assumption is always satisfied if $(x_{t}^a)_{t \geq
0}$ is a stationary on/off process, i.e., if $E=\{0,1\}$. In this
case
%
%\begin{rem}
%If we assume that the state space of the semi-Markov process
%$(x_{t}^a)_{t \geq}$ is binary, i.e., $E=\{0,1\}$, then
\[
    x^{a}_{t}=\frac{\lambda(x^a_{t},s_{t})-\lambda(0,s_{t})}
    {\lambda(1,s_{t})-\lambda(0,s_{t})},
\]
and the representation (\ref{eq:seperable}) holds with
\[
f(x)=x, \quad g(s)=\lambda(1,s)-\lambda(0,s) \quad \text{and}
\quad h(s)=\lambda(0,s).
\]
%\end{rem}
\end{example}

We are now ready to show that the fluctuations of the logarithmic
stock price around its first order approximation behaves like a
fractional Ornstein-Uhlenbeck process.

\begin{thm}
Under the Assumptions~\ref{assumplambda} and
~\ref{ass:separable}  %Assume further that
we have that
\[
    \mathcal{L}\mbox{-}\lim_{T \rightarrow \infty}
    \mathcal{L}\mbox{-}\lim_{N \rightarrow \infty}
    T^{1-H}\frac{\sqrt{N}}{\sqrt{L(T)}} \left\{ S^{N,T}_{t} - s_t \right \}_{0 \leq t \leq 1}
    = \{ \hat{Z}_t \}_{0 \leq t \leq 1}
\]
Here $\hat{Z}$ denotes unique solution to the stochastic
differential equation
\[
    d\hat{Z}_t = \bar{\lambda}'(s_t) \hat{Z}_t dt + \sigma g(s_t) d B^H_t
\]
where $B^H$ is a fractional Brownian motion with Hurst coefficient
$H=\frac{3-\alpha}{2}$. The integral with respect to $B^H$ is
understood as a limit in probability of Stieltjes sums.
\end{thm}
\begin{Proof} The proof uses modifications of arguments given
in the proof of Theorem~\ref{thm-limit-N} and the approximation
result for integrals with respect to fractional Brownian motion in
\cite{kn:bhs}.

\begin{itemize}
\item[(i)] In a first step we study the dynamics of the process
$\{Y^{N,T}_t\}_{t \geq 0}$ defined by
\[
    Y^{N,T}_t = \sum_{a=1}^N \int_0^t \frac{ \lambda(x^a_{Tu},s_u) -
    \bar{\lambda}(s_u)}{\sqrt{N}} \ du
\]
Under Assumption~\ref{ass:separable} we can write
\begin{equation}
\begin{split}
    Y^{N,T}_t %&= \sum_{a=1}^N \int_0^t \frac{ \lambda(x^a_{Tu},s_u) -
%    \bar{\lambda}(s_u)}{\sqrt{N}} \ du \\
    &= \sum_{a=1}^N \int_0^t \frac{1}{\sqrt{N}}
    \left[ f(x^a_{Tu}) g(s_u) + h(s_u) - \bar{\lambda}(s_u) \right] du \\
    &=\sum_{a=1}^N \int_0^t \frac{1}{\sqrt{N}}
    \left[ f(x^a_{Tu}) - \hat{\mu} \right] g(s_u) du
\end{split}
\end{equation}
Since $f$ is one-to-one, $(f(x^{a}_t))_{t \geq 0}$ is a
semi-Markov process that has the same sojourn time structure as
the underlying semi-Markov process $(x^{a}_{t})_{t \geq 0}$. In
particular, $f(0)$ is the state whose sojourn time distribution
has heavy tails. Therefore it follows from Theorem 4.1 of
\cite{kn:bhs} that
\begin{equation}\label{eq:first-step}
    \mathcal{L}\mbox{-}\lim_{T \rightarrow \infty}
    \mathcal{L}\mbox{-}\lim_{N \rightarrow
    \infty} T^{1-H}\left\{ \frac{1}{\sqrt{L(T)}}Y^{N,T}
    %\frac{T^{1-H}}{L(T)} \sum_{a=1}^N \int_0^t
    %\frac{1}{\sqrt{N}} \left( \left(f(x^a_{Tu})-\mu\right)\right)du
    \right\}_{0 \leq t \leq 1} =
    \left\{ \sigma \int_0^t g(s_u) d B^{H}_{u} \right\}_{0 \leq t \leq
    1}
\end{equation}
for some $\sigma > 0$ because $\hat{\mu} \neq f(0)$.

\item[(ii)] Let us now define a family of stochastic processes
$\widetilde{Q}^{N,T} = \{\widetilde{Q}^{N,T}_t\}_{0 \leq t \leq
1}$ by
\[
    \widetilde{Q}^{N,T}_t \triangleq \int_0^t \bar{\lambda}'(s_u)
    \widetilde{Q}^{N,T}_u du + \frac{T^{1-H}}{\sqrt{L(T)}} Y^{N,T}_t
    + \frac{T^{1/2-H}}{\sqrt{L(T)}} X^N_t.
\]
Since the rate functions are bounded and $H > \frac{1}{2}$
\[
    \lim_{T \rightarrow \infty} \sup_N \sup_{0 \leq t \leq 1}
    \frac{T^{1/2-H}}{ \sqrt{L(T)}} X^N_t = 0 %\qquad \Prob\mbox{-a.s.}
\]
almost surely, and the continuous mapping theorem along with (i)
yields
\[
    \mathcal{L}\mbox{-}\lim_{T \rightarrow \infty}
    \mathcal{L}\mbox{-}\lim_{N \rightarrow \infty}
    \{ \widetilde{Q}^{N,T}_t \}_{0 \leq t \leq 1}
    = \{ \hat{Z}_t \}_{0 \leq t \leq 1}
\]

\item[(iii)] Let us put
\[
    Q^{N,T}_t \triangleq T^{1-H} \frac{\sqrt{N}}{\sqrt{L(T)}} (S^{N,T}_t -
    s_t).
\]
Up to a term of the order $\frac{\log N}{\sqrt{N}}$ we obtain
\[%\begin{equation} \label{representation-QN}
    Q^{N,T}_t = \frac{\int_0^t \sum_{a=1}^N
    \left\{ \lambda\left(x^a_{Tu}, S^{N,T}_u\right) - \lambda(x^a_{Tu},s_u)\right\} du}
    {\sqrt{N}} + \frac{T^{1-H}}{\sqrt{L(T)}} Y^N_t + \frac{T^{1/2-H}}{\sqrt{L(T)}} X^N_t.
\]
Using the same arguments as in the proof of Theorem
\ref{thm-limit-N} we thus see that
\[
    \lim_{N \rightarrow \infty} \sup_{0 \leq t \leq 1} |Q^{N,T}_t -
    \widetilde{Q}^{N,T}_t | = 0 \quad \mbox{in probability}
\]
for all $T \in \n$. Hence the assertion follows from (ii).
\end{itemize}
\end{Proof}

\begin{rem}
In the case of Markov switching, i.e., when the process $x_{t}$ is
a Markov process, we obtain standard Ornstein-Uhlenbeck process,
i.e., we have that
\[
    \mathcal{L}\mbox{-}\lim_{T \rightarrow \infty}
    \mathcal{L}\mbox{-}\lim_{N \rightarrow \infty}
    \sqrt{T}\frac{\sqrt{N}}{\sqrt{L(T)}} \left\{ S^{N,T}_{t} - s_t \right \}_{0 \leq t \leq 1}
    = \{ \tilde{Z}_t \}_{0 \leq t \leq 1},
\]
where $\tilde{Z}$ denotes unique solution to the stochastic
differential equation
\[
    d\tilde{Z}_t = \bar{\lambda}'(s_t) \tilde{Z}_t dt + \sigma g(s_t) d
    B_t,
\]
with $B$ a standard Brownian motion.
\end{rem}

%%%%%%%%%%%%%%%%%%%%%%%%%%%%%%%%%%%%%%%%%%%%%%%%%%%%%%%%%%%%%%%%%

\section{Outlook \& Conclusion\label{conc}}
We briefly outline two possible avenues of future research:
microstructure models of fractional volatility and strategic
interactions between ``big players.''

\subsection{Fractional Volatility}

In this article, we suggested a microeconomic approach to
financial price fluctuations that is capable of explaining the
decay of the Hurst coefficient of the S\&P 500 index in the late
1990s. We note that the evidence of long memory in stock price
returns is mixed, there are several papers in the empirical
finance literature providing evidence for the existence of long
memory, yet there are several other papers that contradict these
empirical findings; see e.g. \cite{bps} for an exposition of this
debate and references. However, long memory is a well accepted
feature in volatility (squared and absolute returns) and trading
volume (see e.g. \cite{contstylized} and \cite{engle}). We are now
going to illustrate how the mathematical results of this paper
might also be seen as an intermediate step towards a
microstructural foundation for this phenomenon. To ease notational
complexity and to avoid unnecessary technicalities we restrict
ourselves to the simplest case where the order rates do not depend
on asset prices. Specifically, we assume that (after taking the
$N$-limit) the dynamics of the asset price process can be
described by a stochastic equation of the from
\begin{eqnarray*}\label{eq:mimT1}
    S^{T}_{t} = \frac{1}{T} \left\{
    \Pi_+ \left(T \int_0^{t} \lambda_+(Y^T_u) \,du  \right)
    - \Pi_-\left(T \int_0^{t} \lambda_-(Y^T_u) \,du  \right)
    \right\}
\end{eqnarray*}
where the Gaussian process $Y^T$ defined in
(\ref{def-Gauss-processes}) converges in distribution to a \fbm
process. In view of the strong approximation of Poisson processes
by Browninan motion, and because the rate functions are bounded,
the evolution of prices can be described in terms of an ordinary
differential equation in a random environment generated by a
fractional Brownian motion:
\[
    \mathcal{L}\mbox{-}\lim_{T \rightarrow \infty}
    \{S^T_t \}_{0 \leq t \leq 1} = \{\hat{s}_t\}_{t \leq 0 \leq 1}
    \quad \mbox{where} \quad d \hat{s}_t = \lambda(B^H_t) \,dt.
\]
%where the process $\{\hat{s}\}_{0 \leq t \leq 1}$ is defined in
%terms of the ODE
%\[
%    d \hat{s}_t = \lambda(B^H_t) \,dt.
%\]
The fluctuations around this first order approximation satisfy
\[
    \sqrt{T}\left( S^T_t -  \int_0^t \lambda(Y^T_u)
    \,du \right) =
    B_+ \left( \int_0^t \lambda_+(Y^T_u) \,du \right) -
    B_- \left( \int_0^t \lambda_-(Y^T_u)  \,du
    \right),
\]
up to a term of the order $\frac {\log T}{\sqrt{T}}$. Convergence
of the Gaussian process $Y^T$ to \fbm along with continuity of the
rate functions yields
\[
    \mathcal{L}\mbox{-}\lim_{T \rightarrow \infty}
    \left\{ B_\pm \left( \int_0^t \lambda_\pm(Y^T_u) \,du \right)
    \right\}_{0 \leq t \leq 1} =
    \left\{ \int_0^t \sqrt{\lambda_\pm(B^H_u)} \, dB^\pm_u
    \right\}_{0 \leq t \leq 1}.
\]
Thus, for large $T$, logarithmic asset prices satisfy
\begin{eqnarray*}
    S^T_t & \stackrel{\mathbb{D}}{\approx} & \int_0^t
    \lambda(Y^T_u)\,du +
    \frac{1}{\sqrt{T}} \int_0^t \sqrt{\lambda_+(Y^T_u)} \, dB^+_u -
    \frac{1}{\sqrt{T}} \int_0^t \sqrt{\lambda_-(Y^T_u)} \, dB^-_u \\
    & \stackrel{\mathbb{D}}{\approx} & \int_0^t
    \lambda(B^H_u)\,du +
    \frac{1}{\sqrt{T}} \int_0^t \sqrt{\lambda_+(B^H_u)} \, dB^+_u -
    \frac{1}{\sqrt{T}} \int_0^t \sqrt{\lambda_-(B^H_u)} \, dB^-_u,
\end{eqnarray*}
i.e., the volatility is driven by a \fbm process which is
independent of the Wiener processes $B^+$ and $B^-$. We will
further elaborate on the microstructure of fractional volatility
in a separate paper.

\subsection{Strategic Interactions}
Together with the price taking small investors, it is also
possible to incorporate the effects of large investors who
influence the price. The existence of large agent price effects
has been empirically described in several papers: \cite{kraus},
\cite{hml} and \cite{chan} describe the impacts of institutional
trades on stock prices. In the presence of large agents there is
limited liquidity in the market since the holdings of the stocks
is concentrated in the hands of a few big traders. Trades of ``big
player's'' also affect stock prices due to large order sizes.

\subsubsection{Stochastic Equations in Strategically Controlled
Environments}

Horst \cite{horst-stoch-eqn}, \cite{horst-GEB} provides a
mathematical framework for analyzing linear stochastic difference
equation of the form (\ref{discrete-time-prices}) when the
dynamics of the random environment is simultaneously controlled by
the actions of strategically interacting agents playing a
discounted stochastic game with complete information. In
\cite{horst-stoch-eqn} we considered a simple microstructure
models where small investors choose their current benchmarks in
reaction to the actions taken by some ``big players''. One may,
for example, think of a central bank that tries to keep the ``mood
of the market'' from becoming too optimistic and, if necessary,
warns the market participants of emerging bubbles. One may also
think of financial experts whose recommendations tempt the agents
into buying or selling the stock. These market participants
influence the stock price process through their impact on the
behavior of small investors, but without actively trading the
stock themselves. It seems natural to assume that the big players
anticipate the feedback effect their actions have on the evolution
of stock prices and thus interact in a strategic manner. Under a
weak interaction condition, the resulting stochastic game has a
homogenous Nash equilibrium in Markovian strategies. It turns out
that the main qualitative feature of the models studied in
\cite{Foellmer-Schweizer}, \cite{horstkirman} and \cite{Horst02},
namely namely asymptotic stability of stock prices can be
preserved even in a model of strategic interactions. However, the
long run distribution of stock prices depends on the equilibrium
strategy and is thus not necessarily uniquely determined. Hence,
the presence of strategically interacting market participants can
be an additional source of uncertainty.

\subsubsection{Stochastic Games in a Non-Markovian Setting}

Bayraktar and Poor \cite{bp} considered the strategic interaction
of large investors and found an equilibrium stock price taking
into account that the feedback effects of the large investors on
the price. The large traders find themselves in a random
environment due to the trades of small (i.e. price taking)
investors. In \cite{bp}, the institutional investors strategically
interact through the controls they exert on the coefficients of a
stochastic differential equation driven by a fractional Brownian
motion. Here, the fractional Brownian motion models the effect of
the price taking investors on the price. It can be argued that the
observed stock price is the Nash-equilibrium price that arises as
a result of the strategic interaction of the institutional
investors this random environment. Bayraktar and Poor carries out
an analysis of stochastic differential games in a non-Markov
environment using the stochastic analysis for fractional Brownian
motion developed in \cite{kn:duncan}. This analysis can be viewed
as a first step toward incorporating the feedback effects of the
large investors and the strategic interaction into the description
of the stock price dynamics.

\end{document}